\documentclass[11pt]{article}
\usepackage[T1]{fontenc}
\usepackage[utf8]{inputenc}
\usepackage{lmodern}
\usepackage[margin=1in]{geometry}
\usepackage{amsfonts,amsmath,amssymb,amsthm,mathtools,mathrsfs,times,color}
\usepackage{microtype}
\usepackage{enumitem,booktabs,array}
\usepackage[colorlinks,linkcolor=blue,anchorcolor=blue,citecolor=blue]{hyperref}
\numberwithin{equation}{section}
\newtheorem{theorem}{Theorem}[section]
\newtheorem{lemma}[theorem]{Lemma}
\newtheorem{proposition}[theorem]{Proposition}
\newtheorem{corollary}[theorem]{Corollary}

\theoremstyle{definition}
\newtheorem{remark}[theorem]{Remark}
\newcommand{\ddc}{\sqrt{-1}\partial\bar\partial}
\newcommand{\Ca}{\mathcal C}
\newcommand{\DF}{\operatorname{DF}}
\newcommand{\Aut}{\operatorname{Aut}}
\newcommand{\Vol}{\operatorname{Vol}}

\newcommand{\tr}{\operatorname{tr}}
\newcommand{\reg}{\mathrm{reg}}

\newcommand{\C}{\mathbb C}
\newcommand{\R}{\mathbb R}
\newcommand{\Q}{\mathbb Q}
\newcommand{\Z}{\mathbb Z}
\DeclareMathOperator{\Lie}{Lie}
\setlist[enumerate]{itemsep=4pt,topsep=5pt}

\newcommand{\Stab}{\operatorname{Stab}}
\newcommand{\divisor}{\operatorname{div}}
\newcommand{\Supp}{\operatorname{Supp}}
\newcommand{\orb}{\mathrm{orb}}
\newcommand{\HH}{\mathcal H}
\newcommand{\MM}{\mathcal M}
\newcommand{\D}{\mathcal D}
\DeclareMathOperator{\Impart}{Im}
\hypersetup{pdftitle={Calabi flow on polarized K\"ahler surfaces},
 pdfauthor={Haozhao Li, Bing Wang, Kai Zheng},
 pdfsubject={}}
\title{Calabi flow on polarized K\"ahler surfaces}
\author{Haozhao Li\footnote{Supported by NSFC grant No. 12471058 and a
research fund from University of Science and Technology of China.},
\quad Bing Wang\footnote{Supported by NSFC grant No. 12431003 and a
research fund from University of Science and Technology of China.}
\quad and \quad Kai Zheng\footnote{Supported by NSFC grant No. 12571091.}}
\date{}
\begin{document}
\maketitle
\begin{abstract}
We study global smooth Calabi flows on polarized K\"ahler surfaces
$(M,L)$, assuming uniformly bounded Sobolev constants and the algebraic
conditions $\chi(\mathcal O_M)=1$ and $c_1(M)\cdot c_1(L)>0$.
Using the orbifold compactness theorem, we prove
a uniform Bergman estimate, construct a normal test configuration for
each geometric limit, and establish uniqueness of the extremal orbifold
limit up to polarized holomorphic isometry. K-polystability then
identifies the limit with the original surface and gives a constant
scalar curvature K\"ahler metric. We also prove finite-group equivariant
versions of these results. Independently of the flow, we obtain a
uniform Bergman lower bound for polarized cscK surfaces with bounded
volume and Sobolev constants and $|S|\leq1$.
We illustrate the hypotheses with nonstationary flows near a
product cscK metric on $\mathbb P^1\times\mathbb P^1$ and near
Calabi's non-cscK extremal metric on the first Hirzebruch surface.
Their controlled-cone energy inequalities follow from elementary
calculations and continuity.
\end{abstract}
\begingroup
\small
\setcounter{tocdepth}{2}
\tableofcontents
\endgroup

\section{Introduction}\label{app:sec:intro}
This is the third paper in the series of Li--Wang--Zheng on the Calabi
flow on K\"ahler surfaces. The first paper \cite{LongTime} establishes
an extension criterion in terms of the Sobolev constant. The second
paper \cite{Compactness} proves weak compactness for global flows
with uniformly bounded Sobolev constants.

The purpose of this paper is to establish uniqueness of the geometric
limit of the Calabi flow on polarized K\"ahler surfaces, up to polarized
holomorphic isometry, and to identify this limit using algebraic
stability. Under the analytic and algebraic assumptions stated below,
we prove a uniform Bergman lower bound and use the resulting projective
embeddings to construct test configurations with the geometric limits
as central fibres. K-polystability then identifies the limit with a
constant scalar curvature K\"ahler metric, abbreviated as a cscK metric,
on the original surface.

Let $(M,L)$ be a smooth connected polarized complex surface, and put
$\Omega\coloneqq c_1(L)$. We use the conventions
\begin{equation}\label{app:eq:conventions}
[\omega(t)]=2\pi\Omega,\quad
\partial_t\omega=\ddc S,\quad
\sqrt{-1}F_{h(t)}=\omega(t),\quad dV_t\coloneqq \omega(t)^2/2,
\end{equation}
where $S$ is complex scalar curvature, equal to one half of the
Riemannian scalar curvature. Write
\[
V\coloneqq\int_MdV_t=2\pi^2\Omega^2,\qquad
\bar S\coloneqq V^{-1}\int_M S\,dV_t
       =\frac{2c_1(M)\cdot\Omega}{\Omega^2}.
\]
The volume and average scalar curvature are constant.
The Calabi energy and its volume-normalized form are
\[
\Ca(\omega)\coloneqq\int_M(S_\omega-\bar S)^2dV_\omega,
\qquad \widehat{\Ca}(\omega)\coloneqq\frac{\Ca(\omega)}{V}.
\]
Along the flow we write $\Ca(t)\coloneqq\Ca(\omega(t))$.

We assume a global smooth Calabi flow satisfying two conditions:
the \emph{uniform Sobolev bound}
\begin{equation}\label{tc:eq:sob}
\sup_{t\geq0}C_S(\omega(t))\leq C_S<\infty.
\end{equation}
Here $C_S(\omega(t))$ denotes the least constant such that
\[
\|u\|_{4,t}^2\leq C_S(\omega(t))
\left(\|\nabla u\|_{2,t}^2+V^{-1/2}\|u\|_{2,t}^2\right)
\]
for all $u\in W^{1,2}(M,g(t))$, where the norms use $dV_t$.
The second is the \emph{numerical condition} on the polarized surface:
\begin{equation}\label{tc:eq:numerical}
\chi(\mathcal O_M)=1,\qquad c_1(M)\cdot c_1(L)>0.
\end{equation}

Both conditions are natural. Condition~\eqref{tc:eq:sob} is the
analytic input: a bounded Sobolev constant permits extension across
a finite endpoint \cite[Theorem 1.2]{LongTime} and, for a global flow,
gives subsequential convergence to extremal K\"ahler orbifolds with
isolated quotient singularities \cite{Compactness}. We recall this
compactness property in Section~\ref{app:sec:preliminaries}; the limiting
complex structure may differ from the original one.

Condition~\eqref{tc:eq:numerical} is the algebraic input. Positive
anticanonical degree is equivalent to positive average scalar curvature;
its preservation in the limit gives $H^2(Z,\mathcal O_Z)=0$.
Comparison of holomorphic section spaces, together with
$\chi(\mathcal O_M)=1$, then gives $\chi(\mathcal O_Z)=1$ and
$H^1(Z,\mathcal O_Z)=0$. These identities fix the Hilbert polynomial
and allow us to recover complete limiting section spaces, following
Zhang--Zhou \cite[Section 4]{ZZ}.

Both conditions fit the controlled-cone applications of
\cite[Corollary 1.3]{LongTime}. For a K\"ahler class $\Omega$, the
cohomological inequalities there are
\[
c_1(M)\cdot\Omega>0,\qquad
c_1(M)^2-\frac23\frac{(c_1(M)\cdot\Omega)^2}{\Omega^2}>0.
\]
These inequalities force $M$ to be rational: $K_M\cdot\Omega<0$
makes all plurigenera vanish, so $M$ is rational or ruled by the
Enriques--Kodaira classification; $c_1(M)^2>0$ excludes ruled minimal
models over curves of genus $g\geq1$, whose first Chern number after
$r$ blow-ups is $8(1-g)-r\leq0$. Hence $q(M)=p_g(M)=0$ and
$\chi(\mathcal O_M)=1$, giving \eqref{tc:eq:numerical} whenever
$\Omega$ lies on the ray of $c_1(L)$. The additional initial-energy
inequality in \cite[Corollary 1.3]{LongTime} gives global existence and
\eqref{tc:eq:sob}; it must still be checked for the chosen metric.

Global smooth Calabi flows on complex surfaces are also known under
several other assumptions. Chen--He \cite[Theorem 1.2]{CHFlow} and
Huang--Zheng \cite[Theorem 1.4]{HZFlow} establish stability results
near cscK and extremal metrics, respectively. Tosatti--Weinkove
\cite[Theorems 1--2]{TWFlow} obtain global existence from small initial
Calabi energy and initial curvature bounds when $c_1=0$, and when
$c_1<0$ for classes close to the canonical class. Under symmetry
assumptions, Chen--He \cite[Theorem 1.1]{CHToric} treat toric Fano
surfaces with positive extremal Hamiltonian potential and an
initial-energy bound; Guan \cite[Theorem A and Theorem 10.2]{GuanFlow}
treats certain $\mathbb P^1$-bundles, proving convergence of the
modified flow; and Feng--Huang \cite[Theorems 1.2--1.3]{FHFlow}
prove global existence and convergence for $T^2$-invariant metrics
in the standard flat class on
$\mathbb C^2/(\mathbb Z^2+i\mathbb Z^2)$.

Sz\'ekelyhidi's examples \cite[Theorems~1--2 and p.~839]{SzekRuled}
illustrate degeneration at infinity along smooth global Calabi flows.
For momentum
metrics on a ruled surface $X$ over a genus-two curve, he proves
that, in classes admitting no extremal metric, the flow degenerates
as $t\to\infty$ into complete extremal pieces, with possible collapse.
These flows have unbounded Sobolev constants and
$\chi(\mathcal O_X)=-1$. In his notation $L_m\coloneqq C+mS_\infty$, with
$C$ a fibre and $S_\infty$ a section with $S_\infty^2=1$, adjunction gives
$c_1(X)\cdot L_m=2-m<0$ for $m\geq k_1\simeq18.9$.
Thus both \eqref{tc:eq:sob} and \eqref{tc:eq:numerical} fail.

To state the main result, for an $L^2(dV_t,h(t)^k)$-orthonormal basis
of $H^0(M,L^k)$, define
\begin{equation}\label{pc:eq:bergman}
\rho_{k,t}(x)\coloneqq \sum_a|s_a(x)|_{h(t)^k}^2
=\sup_{s\ne0}\frac{|s(x)|_{h(t)^k}^2}{\|s\|_{2,t}^2}.
\end{equation}
Here K-semistability means that every normal relatively ample test
configuration for $(M,L)$, of any positive exponent, has nonnegative
Donaldson--Futaki invariant, denoted by $\DF$; see
\cite[Definition 3.11 and Corollary 3.16]{BHJ}.
K-polystability additionally requires every zero-invariant configuration
to be a product, allowing a nontrivial $\C^*$-action on $M$.
We use the unmodified Donaldson--Futaki invariant, with no relative
correction or uniform positive lower bound assumed.

\begin{theorem}[Applications of compactness]\label{app:thm:main}
Suppose a global smooth Calabi flow $\omega(t)$ on $(M,L)$ satisfies the
Sobolev bound \eqref{tc:eq:sob} and the numerical condition
\eqref{tc:eq:numerical}. There is a common integer $r>0$ with the
following properties. Geometric limits below are subsequential limits
as $t\to\infty$ in the sense of Theorem~\ref{app:thm:compactness},
evaluated at translated time zero.
\begin{enumerate}[font=\normalfont\bfseries]
\item \textbf{\upshape Bergman estimate.} There are $0<b\leq B<\infty$ such that
\[
b\leq\rho_{r,t}(x)\leq B\qquad(t\geq0,\ x\in M).
\]
On every orbifold limit, the limiting orbifold line bundle $L_Z$ has
$L_Z^r$ descending to a very ample Cartier line bundle $A_Z$ on the
coarse surface $Z$. Along any geometrically
convergent subsequence, after further extraction, complete orthonormal
bases at this power can be chosen to converge smoothly on regular
compact subsets to a complete orthonormal basis on the limit.
\item \textbf{\upshape Test configurations.} Every polarized geometric limit
$(Z,A_Z)$ is the central polarized fibre of a normal, relatively ample
test configuration for $(M,L)$ of exponent $r$. One can choose this
configuration with strictly negative Donaldson--Futaki invariant when
the limiting metric is not cscK, and with zero invariant when it is cscK.
\item \textbf{\upshape Uniqueness of the Calabi flow limit.}
This fixed Calabi flow has a unique polarized extremal orbifold limit
$(Z_\infty,A_\infty,\omega_\infty)$ up to polarized holomorphic isometry,
with $[\omega_\infty]=(2\pi/r)c_1(A_\infty)$.
Explicitly, for any two geometric limits $(Z_i,A_i,\omega_i)$,
$i=1,2$, there is an orbifold biholomorphism $F:Z_1\to Z_2$ satisfying
\begin{equation}\label{un:eq:isom}
F^*A_2\simeq A_1,\qquad F^*\omega_2=\omega_1.
\end{equation}
The Calabi flow converges to this isometry class in the quotient
geometric topology of Theorem~\ref{app:thm:compactness}. In particular,
$(M,g(t))$ converges to $(Z_\infty,g_\infty)$ in the Gromov--Hausdorff
(GH) sense. Every sequence $t_j\to\infty$ has a subsequence for which the
metrics and complex structures converge smoothly on regular compact
subsets of $Z_\infty$, after pullback by embeddings into $M$.
These embeddings may depend on the subsequence.
\item \textbf{\upshape Stability and existence.}
If $(M,L)$ is K-semistable, then $S_{\omega_\infty}\equiv\bar S$;
the cscK orbifold $Z_\infty$ may still be singular.
If $(M,L)$ is K-polystable, then $Z_\infty$ is smooth and there is
a polarized isomorphism
\[
F:(M,L^r)\longrightarrow(Z_\infty,A_\infty).
\]
The pullback $\omega_{\mathrm{cscK}}\coloneqq F^*\omega_\infty$ is a
smooth cscK metric on the original complex surface $M$, with
$[\omega_{\mathrm{cscK}}]=2\pi c_1(L)$. The Calabi flow converges to
the polarized holomorphic-isometry class of $(M,L^r,\omega_{\mathrm{cscK}})$
in the sense of part~(3).
\end{enumerate}
\end{theorem}

The finite-group equivariant version, Theorem~\ref{eqv:thm:main},
chooses the limiting embeddings and test configurations equivariantly
when the group acts on the polarization and preserves the initial
metric. The symmetry is preserved by the flow and carried through
the compactness and finite-dimensional constructions. Since the
resulting central fibres are normal integral surfaces, to obtain a
smooth cscK metric it suffices to check
nonnegativity of the Donaldson--Futaki invariant on equivariant normal
ample test configurations with integral central fibre, with equality
only for products; see Corollary~\ref{eqv:cor:restricted}. This restricted
criterion makes the equivariant theorem useful when symmetry gives
access to stability inequalities for these configurations.

Section~\ref{ex:sec:examples} gives two examples of nonstationary
flows satisfying the hypotheses. On
$(\mathbb P^1\times\mathbb P^1,\mathcal O(1,2))$, a small potential
perturbation of the product cscK metric converges to its polarized
isometry class. On $(\mathbb F_1,-K_{\mathbb F_1})$, a small
$U(2)$-invariant perturbation of Calabi's extremal metric converges,
after holomorphic pullback, to a non-cscK extremal metric. The
controlled-cone inequalities and the extremal energy are computed
exactly; neither example requires numerical integration. Their
existence and local flow convergence are known from
\cite{CHFlow,HZFlow}; they illustrate the conclusions of our theorem.

We outline the proof of Theorem~\ref{app:thm:main}. To compare the
geometric limits algebraically, we first embed them by complete section
spaces of a common degree. The limiting polarization extends as an
orbifold line bundle $L_Z$. Section upper bounds, regular convergence
and uniform volume tails preserve section inner products, initially
giving an injection into the limiting section space.

The Sobolev bound controls the local quotient-group orders, so one
integer $q$ makes every $L_Z^q$ descend to an ample Cartier bundle
$B_Z$. The numerical argument following Zhang--Zhou
\cite[Corollaries 2.11--2.12 and Remark 4.7]{ZZ} gives
\[
P_q(k)\coloneqq\chi(M,L^{qk})=\chi(Z,B_Z^k).
\]
Matsusaka's theorem for normal varieties with rational singularities
\cite[Theorem 2.4, pp.~206--207; see also p.~176]{Matsusaka} supplies
a common very ample power $m_0$, with vanishing higher cohomology,
depending only on this polynomial. Thus the power $r\coloneqq qm_0$
is fixed before choosing a subsequence; on each limit set
$A_Z\coloneqq B_Z^{m_0}$.
Equality of section dimensions makes the limiting orthonormal bases
complete. Their convergence on the regular locus embeds the geometric
limit as a closed subscheme of the Hilbert limit; equality of their
full Hilbert polynomials excludes any additional scheme structure.

The Bergman lower bound uses the isolated nature of the exceptional
points in complex dimension two. We choose a smooth ample divisor on
the limit avoiding those finitely many points, and approximate its
defining section by sections along the flow. Connectedness of the
approximating zero divisors excludes additional zeros in the exceptional
regions. The minimum principle for the logarithm of the section norm
then gives a positive lower bound there, while complete-section
convergence controls the remaining regular region. Because the same
power $r$ works for every subsequence, a contradiction argument makes
this bound uniform for all large times; positivity and continuity
handle bounded time intervals. Since the argument does not use
H\"ormander's $L^2$ estimate, no lower bound for the Ricci curvature
is needed.

The Bergman argument uses only polarized weak compactness, convergence
of complete section spaces and uniform section upper bounds.
It therefore also applies to families of polarized cscK surfaces
with bounded volume and Sobolev constants and $|S|\leq1$.
This answers, in complex dimension two under these bounds, the
long-standing question concerning a uniform Bergman lower bound for
cscK metrics; see Zhang--Zhou \cite[Problem 8.1]{ZZ} for its explicit
formulation and Remark~\ref{pc:rem:zz} for further details.

Next, fix an ordered basis $s\coloneqq (s_1,\ldots,s_N)^T$ of $H^0(M,L^r)$,
where $N\coloneqq h^0(M,L^r)$. Its Gram matrix with respect to the natural
$L^2$ inner product is the Hermitian matrix $G(t)$ with entries
\[
G(t)_{ab}\coloneqq H_{r,t}(s_a,s_b)
\coloneqq \int_M h(t)^r(s_a,s_b)\,dV_t,
\qquad 1\leq a,b\leq N.
\]
Hermitian pairings are linear in their first entry. On bounded
translated time intervals, normalized Gram matrices converge to
$e^{-2\tau C_Z}$, where $C_Z$ represents the extremal flow.
Futaki--Mabuchi rationality \cite{FM} and connectedness of the spectral
cluster set give a common trace-free spectrum, represented by
$\Lambda$. Chen--Sun--Wang's finite-dimensional argument
\cite[Section 2]{CSW} then produces a compact connected Hilbert
limit set $\mathscr C$ in the closure of one orbit of
$G_\Lambda\coloneqq\operatorname{Cent}_{GL(N,\C)}(\Lambda)$.

Rationality makes the first degeneration algebraic. Reductivity of
the stabilizer of a prescribed metric limit supplies a second
degeneration in $G_\Lambda$. The two subgroups commute, so their
product, with a sufficiently large power of the first, gives one
test configuration whose central fibre is the polarized variety
underlying that geometric limit. Its Donaldson--Futaki invariant
is negative for a non-cscK limit and zero for a cscK limit.
This proves the stability consequences using complete-section
convergence, independently of metric-limit uniqueness or the global
Bergman lower bound.

For uniqueness, put $K_\Lambda\coloneqq G_\Lambda\cap U(N)$.
Fixed-orbifold extremal uniqueness, proved in
Subsection~\ref{app:sec:orbifold-uniqueness}, gives
$\mathscr C\cap G_\Lambda[Z]=K_\Lambda[Z]$.
Together with reductivity and connectedness, this forces
$\mathscr C$ to be one $K_\Lambda$-orbit. Complete orthonormal
bases then identify the polarized metrics as well as the varieties;
translated-window convergence extends the conclusion to arbitrary
time sequences.

Constants and gauges may depend on the fixed flow. Its metric-limit
uniqueness requires no stability and makes no assertion about different
initial metrics, bubble trees, or uniqueness of test configurations.
The combined configuration need not induce the same weighted filtration
as the Calabi flow, and we assert no optimality property; see
Remark~\ref{tc:rem:two-stage}.

The paper is organized as follows.
Section~\ref{app:sec:preliminaries} supplies the compactness and
fixed-orbifold uniqueness results. Section~\ref{app:sec:bergman}
treats complete section spaces and Bergman estimates, including the
uniform bound for families of polarized cscK surfaces.
Section~\ref{app:sec:convergence} develops the finite-dimensional
argument, test configurations, limit uniqueness and stability consequences.
Section~\ref{eqv:sec:main} gives the equivariant results, and
Section~\ref{ex:sec:examples} presents the two smooth-limit examples.

\medskip
\noindent\textbf{On the usage of AI.}
The authors used ChatGPT as an assisting tool in the development of
the paper. The authors take full responsibility for the contents of
the paper.

\par\addvspace{\baselineskip}
\noindent\textbf{Acknowledgement.}
The authors are deeply grateful to Professor Xiuxiong Chen
for his encouragement and constant support.
Bing Wang would like to thank Professor Song Sun for his interest in this work.

\section{Preliminaries}\label{app:sec:preliminaries}

We first collect the principal notation, symbols and terminology.
We then recall the compactness theorem and its equivariant refinement,
and establish the fixed-orbifold extremal uniqueness used to identify
the projective limits, following the convex perturbation method of
Berman--Berndtsson.

\subsection{Notation, terminology and constants}
\label{app:index:notation}
For the convenience of readers, we collect the principal notation,
symbols and terminology in the following table. Local notation is
specified where it occurs.

\noindent\emph{Flow and geometric limits.}
\begin{description}[style=nextline,leftmargin=2.7cm,labelwidth=2.4cm,labelsep=0.3cm,align=left,font=\normalfont,itemsep=0pt,parsep=0pt,topsep=4pt]
\item[$(M,L)$] The smooth polarized complex surface.
\item[$\Omega$] The class $c_1(L)=[\omega(t)]/(2\pi)$.
\item[$\omega(t),g(t)$] The K\"ahler form and its Riemannian metric,
with $\partial_t\omega=\ddc S$.
\item[$h(t)$] The Hermitian metric on $L$ with
$\sqrt{-1}F_{h(t)}=\omega(t)$.
\item[$S,\bar S$] Complex scalar curvature and its spatial average;
$S$ is one half of the Riemannian scalar curvature.
\item[$dV_t,V$] The volume form $\omega(t)^2/2$ and the fixed total
volume $V=2\pi^2\Omega^2$.
\item[$\Ca,\widehat{\Ca}$] The Calabi energy
$\Ca(\omega)=\int_M(S_\omega-\bar S)^2\,dV_\omega$ and its
volume-normalized form $\widehat{\Ca}=\Ca/V$;
$\Ca(t)=\Ca(\omega(t))$.
\item[$C_S$] The uniform scale-invariant Sobolev bound
\eqref{tc:eq:sob}.
\item[cscK, GH] Constant scalar curvature K\"ahler, and
Gromov--Hausdorff.
\item[Compactness property] The time-independent GH comparisons,
regular smooth convergence and volume convergence in
Theorem~\ref{app:thm:compactness}.
\item[$Z_{\orb},Z,\pi$] The limiting orbifold, its coarse complex space,
and the natural map $\pi:Z_{\orb}\to Z$; see
Subsection~\ref{app:sec:geometric-compactness}.
\item[$Z_{\reg}$] The regular locus of the limiting space.
\item[$W,F_\tau$] The real extremal field
$W=\frac12\nabla_{\mathrm{real}}S$ and its flow,
with $\omega_Z(\tau)=F_\tau^*\omega_Z(0)$.
\end{description}

\par\addvspace{\baselineskip}
\noindent\emph{Polarization and section spaces.}
\begin{description}[style=nextline,leftmargin=2.7cm,labelwidth=2.4cm,labelsep=0.3cm,align=left,font=\normalfont,itemsep=0pt,parsep=0pt,topsep=4pt]
\item[$L_Z,h_Z$] The limiting orbifold line bundle and its positive
Hermitian metric.
\item[$B_Z,A_Z$] The ample and very ample Cartier bundles on the coarse
space, with $\pi^*B_Z\simeq L_Z^q$ and $A_Z=B_Z^{m_0}$.
\item[$q,m_0,r$] The common descent power, very ample power and
embedding exponent, with $r=qm_0$; see
Subsection~\ref{app:sec:sections}.
\item[$P_q,P$] The fixed Hilbert polynomials
$P_q(k)=\chi(M,L^{qk})=\chi(Z,B_Z^k)$ and $P(k)=P_q(m_0k)$.
\item[$H_{r,t},N$] The section inner product
$H_{r,t}(s,s')=\int_Mh(t)^r(s,s')\,dV_t$ and
$N=h^0(M,L^r)$; pairings are linear in the first entry.
\item[$\rho_{k,t}$] The Bergman function, the sum of squared pointwise
norms of a complete orthonormal section basis; see
\eqref{pc:eq:bergman}.
\end{description}

\par\addvspace{\baselineskip}
\noindent\emph{Finite-dimensional dynamics and stability.}
\begin{description}[style=nextline,leftmargin=2.7cm,labelwidth=2.4cm,labelsep=0.3cm,align=left,font=\normalfont,itemsep=0pt,parsep=0pt,topsep=4pt]
\item[$G(t)$] The Gram matrix of a fixed section basis for $H_{r,t}$;
this is distinct from an algebraic group denoted by $G$.
\item[$C_Z,\Lambda$] The Hermitian generator with limiting Gram matrix
$e^{-2\tau C_Z}$, and the common representative, up to unitary
conjugation, of its trace-free part.
\item[$C_0$] A fixed change of evaluation coordinates identifying the
flow-induced filtration with the weight filtration of $\Lambda$;
see Lemma~\ref{app:lem:limitset}.
\item[{$[M],[\bar M]$}] The Hilbert point of the chosen complete embedding,
and its intermediate limit
$[\bar M]=\lim_{t\to\infty}e^{t\Lambda}C_0[M]$.
\item[$\mathscr C$] The compact connected compatible Hilbert limit set
of Lemma~\ref{app:lem:limitset}.
\item[$G_\Lambda,K_\Lambda$] The centralizer
$\operatorname{Cent}_{GL(N,\C)}(\Lambda)$ and its compact subgroup
$G_\Lambda\cap U(N)$.
\item[$G,K$] In the abstract orbit lemma, a complex reductive group
and a maximal compact subgroup; see Lemma~\ref{un:lem:orbit}.
\item[$\lambda,\mu,\nu_m$] The commuting algebraic one-parameter
subgroups used to reach the metric limit, with
$\lambda(z)=\exp(-b\Lambda\log z)$ and $\nu_m=\lambda^m\mu$;
see \eqref{tc:eq:lambda}--\eqref{tc:eq:nu}.
\item[$\Aut(Z,A,Y)$] The automorphisms preserving the polarization
$A$ and commuting with the extremal field $Y$; see
Subsection~\ref{app:sec:orbifold-uniqueness}.
\item[$d_k,w_k(\eta)$] The dimension of $H^0(Z,A^k)$ and the total
weight of an action $\eta$ on that section space.
\item[$\DF$] The Donaldson--Futaki invariant, using section weights;
the normalization and sign are fixed in
Subsection~\ref{app:sec:tc}.
\end{description}

\par\addvspace{\baselineskip}
\noindent\emph{Equivariance and examples.}
\begin{description}[style=nextline,leftmargin=2.7cm,labelwidth=2.4cm,labelsep=0.3cm,align=left,font=\normalfont,itemsep=0pt,parsep=0pt,topsep=4pt]
\item[$\Gamma$] A finite local group in an orbifold chart; in
Section~\ref{eqv:sec:main}, the fixed finite global symmetry group
acting on the polarized surface.
\item[$\mathcal R,K_{\mathcal R}$] The centralizer of the global
$\Gamma$-representation on $H^0(M,L^r)^*$ and its unitary subgroup;
see \eqref{eqv:eq:centralizer}.
\item[$\delta(A)$] The controlled-cone quantity
$c_1(X)^2-2(c_1(X)\cdot A)^2/(3A^2)$, for $A=c_1(L)$;
see \eqref{ex:eq:general-threshold}.
\end{description}

\subsection{Geometric compactness}\label{app:sec:geometric-compactness}

We first fix the distinction between an orbifold and its underlying
complex space. The orbifolds occurring here have compatible local
charts $(\widetilde U,\Gamma)$, where $\widetilde U\subset\C^2$ is a
ball and $\Gamma\subset U(2)$ is finite and acts freely away from
the origin. The \emph{coarse space} is the normal complex analytic
surface obtained by gluing the quotients $U\coloneqq \widetilde U/\Gamma$,
with structure sheaf
\[
 \mathcal O_Z|_U\coloneqq (q_*\mathcal O_{\widetilde U})^\Gamma,
 \qquad q:\widetilde U\longrightarrow U.
\]
Thus holomorphic functions on $U$ are precisely the invariant
holomorphic functions upstairs. The coarse space has the same
underlying points and topology as the orbifold, but the chart groups
and their actions are not included as additional structure. In
particular, passing to the coarse space does not resolve its quotient
singularities. When necessary we distinguish the two objects by
$\pi:Z_{\mathrm{orb}}\to Z$; otherwise we write $Z$ for both.
Smooth orbifold metrics and forms are defined by their smooth
invariant lifts to the charts. The structure sheaf $\mathcal O_Z$,
its cohomology and, once projectivity is established, Hilbert points
refer to the coarse space.

\emph{Orbifold line bundles.}
An ordinary line bundle $A$ is locally trivial over the coarse
space itself, including its singular points. An \emph{orbifold line
bundle} $\mathcal L$ consists instead of compatible equivariant
line bundles on the uniformizing covers. Locally its model is
\[
\widetilde U\times\C,\qquad
\gamma\cdot(z,v)=(\gamma z,\chi(\gamma)v),\qquad
\chi:\Gamma\longrightarrow\C^*.
\]
Here $\chi$ is a character describing the action on the fibre over
the origin. Orbifold holomorphic sections over $U$ correspond to
holomorphic functions satisfying
$f(\gamma z)=\chi(\gamma)f(z)$ upstairs.

The precise \emph{descent criterion} is
\[
 \mathcal L\simeq\pi^*A\ \text{for a line bundle }A\text{ on }Z
 \quad\Longleftrightarrow\quad
 \text{every local stabilizer acts trivially on its fibre}.
\]
Necessity follows by pulling back a local frame of $A$. Conversely,
when the fibre action at the centre of a chart is trivial, averaging
a local holomorphic frame over the finite group gives an invariant
frame that is nonzero at the centre and hence on a smaller chart.
The transition functions between these frames are invariant and
descend to the quotient, defining $A$.

For example, on the orbifold chart $\C^2/\{\pm1\}$ consider the
line bundle with action
\[
 (z_1,z_2,v)\longmapsto(-z_1,-z_2,-v).
\]
Its sections satisfy $f(-z)=-f(z)$ and vanish at the origin upstairs.
The nontrivial fibre character prevents descent, whereas its tensor
square descends. Away from orbifold points the two bundle notions
coincide.

On a compact orbifold, a sufficiently divisible tensor power of any
orbifold line bundle has trivial stabilizer actions: finitely many
charts suffice, and one can take a common multiple of the orders of
their fibre characters. This explains the distinction below between
the orbifold polarization $L_Z$ and a Cartier polarization $A_Z$ with
$L_Z^r\simeq\pi^*A_Z$. Since the local actions have no
codimension-one stabilizers, the orbifold canonical bundle corresponds
to the $\Q$-Cartier canonical divisor $K_Z$ of the coarse surface;
their rational intersection numbers agree.

We shall use the following form of the compactness theorem throughout
the paper. Its time-independent convergence maps allow us to compare
the polarization and the section spaces on bounded time intervals.

\begin{theorem}[Geometric input from the compactness paper]\label{app:thm:compactness}\label{tc:hyp:wc}
The long-time compactness theorem of \cite[Theorem 1.1 and Definition 6.1]{Compactness}
has the following consequence. If the Calabi flow is smooth on
$[0,\infty)$ and satisfies \eqref{tc:eq:sob}, then any $t_j\to\infty$
has a subsequence converging on every bounded translated time interval
to an eternal extremal flow on one connected compact K\"ahler orbifold $Z$.
The orbifold has finitely many isolated quotient points, with effective
local groups free on the unit three-sphere. Global GH convergence is
realized by time-independent correspondences whose distortions tend to
zero uniformly on bounded time intervals. Convergence is mixed
$C^\infty$ on regular compact subsets under compatible embeddings
$\Phi_j$, independent of time, whose graphs lie in those correspondences.
These embeddings include collars and cover the corresponding regular
regions. The limiting metrics extend smoothly in orbifold charts and
\[
\omega_Z(\tau)=F_\tau^*\omega_Z(0),\qquad
F_\tau\coloneqq \exp\left(\tfrac\tau2\nabla S_Z\right).
\]
The natural volume measures converge and have total mass $V$.
\end{theorem}
We refer to this conclusion as the compactness property. The cited
theorem applies to any fixed real K\"ahler class. In the present paper,
the polarization and the numerical assumptions \eqref{tc:eq:numerical}
are used to derive its algebraic and metric consequences.

The following consequence of volume convergence will be used when
passing scalar-curvature integrals and section inner products to a limit.

\begin{lemma}[Uniform volume tails]\label{app:lem:tails}
For an exhaustion $K_\ell\Subset Z_{\reg}$ by smooth compact domains and
any bounded interval $I$, the compactness property gives
\begin{equation}\label{tc:eq:volume}
\lim_{\ell\to\infty}\limsup_{j\to\infty}\sup_{\tau\in I}
\Vol_{g(t_j+\tau)}(M\setminus\Phi_j(K_\ell))=0.
\end{equation}
\end{lemma}
\begin{proof}
The metrics $\omega_Z(\tau)$ form a smooth family on the compact orbifold
for $\tau\in I$. They are uniformly comparable to a fixed orbifold metric,
so $\sup_{\tau\in I}\Vol_{\omega_Z(\tau)}(Z\setminus K_\ell)\to0$.
For a fixed $\ell$, mixed smooth convergence gives uniform convergence of
the integrals on $K_\ell$. Subtract these integrals from the common total
volume $V$, first let $j\to\infty$, and then let $\ell\to\infty$.
\end{proof}

The remaining observation fixes the scalar-curvature normalization
on each limiting extremal orbifold.

\begin{lemma}\label{tc:lem:extremal}
Every limit supplied by the compactness property is extremal and
$\omega_Z(\tau)=F_\tau^*\omega_Z(0)$, where $F_\tau$ is a complete
one-parameter group of orbifold biholomorphisms. Moreover the scalar
average of the limit is $\bar S$.
\end{lemma}
\begin{proof}
Extremality and the self-similar formula are conclusions of
Theorem~\ref{app:thm:compactness}. The generator is
$W\coloneqq \frac12\nabla_{\mathrm{real}}S$, with $\mathcal L_W\omega=\ddc S$.
We verify the scalar normalization, which will be used in Riemann--Roch.

Along the Calabi flow, we have
\[
\Ca'(t)=-2\int_M|\bar\partial\nabla^{1,0}S|^2dV_t\leq0.
\]
It follows that the scalar $L^2$ norms of the original flow are uniformly
bounded, since $\int_M S^2dV_t=\Ca(t)+\bar S^2V$. Cauchy--Schwarz and
\eqref{tc:eq:volume} show that the scalar integrals over the omitted sets
tend to zero. Passing to the limit on regular compact subsets preserves
$\int S\,dV$ and total volume, and hence preserves $\bar S$.
\end{proof}

We shall also use the finite-group equivariant refinement of the
compactness theorem, whose geometric proof is given in the companion
paper.

\begin{lemma}[Finite-group equivariant compactness, Theorem 6.12 of \cite{Compactness}]\label{eqv:lem:compactness}
Suppose $(M_j,g_j(\tau),J_j)$ converges with the coherent GH and
regular smooth convergence of Theorem~\ref{app:thm:compactness} on
every bounded time interval. If a fixed finite group $\Gamma$ acts
by time-independent biholomorphisms preserving every $g_j(\tau)$,
then, after a further subsequence, its actions converge to one
time-independent action on the limiting orbifold $Z$ by
biholomorphisms preserving every limiting metric. The GH
correspondences can be chosen $\Gamma$-equivariantly, and the local
action maps converge smoothly on regular compact subsets.
\end{lemma}

\subsection{Extremal uniqueness on a fixed orbifold}\label{app:sec:orbifold-uniqueness}

We establish the fixed-orbifold version of extremal uniqueness used
in Subsection~\ref{app:sec:unique}. On smooth compact K\"ahler
manifolds, extremal metrics in a fixed K\"ahler class are unique
modulo holomorphic automorphisms. For the uniqueness theorem and
its analytic foundations, we refer to the work of Chen--Tian
\cite{CT08}, Ross--Witt Nystr\"om \cite{RWN},
Chen--P\u{a}un--Zeng \cite{CPZ}, and Berman--Berndtsson \cite{BB}.

Here we follow the weak-geodesic convexity and
perturbation argument of Berman--Berndtsson
\cite[Propositions 4.3, 4.7, 4.13--4.14 and the proof of Theorem 4.16]{BB}.
We first establish the weak-geodesic and convexity results on local
uniformizing covers, then compare two extremal metrics with the same
extremal field. The orbifold remains fixed throughout the argument.

For a compact connected polarized orbifold $(Z,A)$, let $\Aut(Z)$
denote its group of orbifold biholomorphisms. Its \emph{polarized
automorphism group} is
\[
\Aut(Z,A)\coloneqq\{f\in\Aut(Z):f^*A\simeq A\},
\]
where the isomorphism is one of holomorphic line bundles (orbifold line
bundles when $A$ is an orbifold polarization). 
For a holomorphic vector field $Y$, we write
\[
\Aut(Z,A,Y)\coloneqq\{f\in\Aut(Z,A):f_*Y=Y\},
\]
and $\Aut^0(Z,A,Y)$ denotes the identity component of this centralizer.
When $A$ is a very ample line bundle on the coarse space, the complete
linear system identifies $\Aut(Z,A)$ with the projective stabilizer
of $Z\subset\mathbf P(H^0(Z,A)^*)$. 

\begin{theorem}[The centralizer form of extremal uniqueness]\label{un:thm:fixed}
Let $(Z,A)$ be a compact connected polarized K\"ahler orbifold with
isolated quotient points and $q(Z)=0$. Let $\omega_0,\omega_1$ be two
orbifold-smooth extremal metrics in the same positive multiple of
$c_1(A)$. Suppose that their holomorphic extremal fields coincide:
\[
Y\coloneqq \nabla^{1,0}_{\omega_0}S_0
 =\nabla^{1,0}_{\omega_1}S_1.
\]
Then some $a\in\Aut^0(Z,A,Y)$ satisfies $a^*\omega_1=\omega_0$.
\end{theorem}

\subsubsection{Analytic tools}\label{app:sec:analytic}

We adapt weak-geodesic existence, convexity and the fixed-measure
estimate to local uniformizing covers of a fixed orbifold.

Let $Z$ be a fixed compact connected effective K\"ahler
orbifold with isolated quotient points, and let $\omega$ be a fixed
orbifold-smooth K\"ahler form. Write $n\coloneqq \dim_{\C}Z$. All smoothness and differential estimates
mean smoothness or estimates on local uniformizing covers. Set
\[
\omega_u\coloneqq \omega+\ddc u,\qquad
\HH\coloneqq \{u\in C^\infty_\orb(Z,\R):\omega_u>0\},\qquad
\mathsf V\coloneqq \int_Z\omega^n.
\]
Throughout Subsection~\ref{app:sec:orbifold-uniqueness}, we use
$\omega_u^n$, without the factor $1/n!$.
Orbifold Hodge theory gives the $\partial\bar\partial$ lemma, so two
cohomologous orbifold K\"ahler forms have potentials in this space.
Here and below integration over a chart $\widetilde U/\Gamma$ means
$|\Gamma|^{-1}$ times integration on the cover, with an invariant
partition of unity. Thus Stokes' formula has no boundary contribution
at the quotient points.

\begin{lemma}[Weak geodesics on the fixed orbifold]\label{un:lem:geod}
Two points of $\HH$ are joined by a unique weak geodesic $u_t$ whose
complex Hessian on $Z\times([0,1]\times S^1)$ is bounded in orbifold
charts. In particular, for a fixed product background $\Omega$,
\begin{equation}\label{un:eq:geod-bounds}
0\leq \pi_Z^*\omega+\ddc U\leq C\Omega,\qquad
\|U\|_{C^1}\leq C,\qquad
(\pi_Z^*\omega+\ddc U)^{n+1}=0.
\end{equation}
These bounds can be chosen uniformly for endpoints in a smooth compact
family of strictly positive metrics. If both endpoints are invariant
under a compact group preserving $\omega$, the geodesic is invariant.
\end{lemma}

\begin{proof}
We adapt the Dirichlet argument to the orbifold. Write
$w\coloneqq t+\sqrt{-1}s$, with $s$ periodic. The function
\[
\underline U\coloneqq (1-t)u_0+tu_1+B(t^2-t)
\]
is a strict subsolution for $B$ sufficiently large: its horizontal
block is $(1-t)\omega_{u_0}+t\omega_{u_1}>0$, and its positive
time-time block dominates the bounded mixed blocks. Put
$\underline\Omega\coloneqq \pi_Z^*\omega+\ddc\underline U>0$ and solve
\begin{equation}\label{un:eq:eps}
(\pi_Z^*\omega+\ddc U_\epsilon)^{n+1}
=\epsilon\,\underline\Omega^{n+1},\qquad
U_\epsilon|_{t=0,1}=u_0,u_1,\quad 0<\epsilon\leq1.
\end{equation}
At $\epsilon=1$ the displayed subsolution is a solution. We apply the
continuity method in the orbifold H\"older spaces. The
linearized Dirichlet operator is a positive elliptic Laplacian with
zero boundary data, hence is invertible by the maximum principle and
the Fredholm alternative. The latter follows by taking local elliptic
parametrices on finite covers, averaging them over the local finite
groups, and patching by an invariant partition of unity.

The Dirichlet estimates of \cite[Theorems 1.3--1.4 and Section 3]{Blocki},
following \cite{Chen}, give uniform gradient and complex Laplacian
bounds. Their constants are independent of $\epsilon$: the density is
$f_\epsilon=\epsilon f_1$, so the required derivatives of
$f_\epsilon^{1/(n+1)}$ and the lower bound for
$f_\epsilon^{1/n}\Delta\log f_\epsilon$ are uniform. The strict
subsolution controls the boundary mixed derivatives, and the fixed
positive tangential block controls the normal-normal second derivative through the
determinant equation.

These local estimates apply on the finite uniformizing covers.
A maximum at a quotient point lifts to an ordinary maximum; at the
time boundary the group acts only on the $Z$ variables, leaving the
lifted boundary flat. Finite charts therefore give the same global
bounds. For fixed $\epsilon>0$, the determinant and trace bounds give
uniform ellipticity, and the usual local regularity closes the
continuity method in orbifold H\"older spaces.

The uniform trace bound gives a bound for the full background
Laplacian. Linear interior and boundary $W^{2,p}$ estimates, with the
smooth prescribed boundary values, give compactness in $C^{1,\alpha}$
for $\alpha<1$. A subsequence tends to $U$ satisfying
\eqref{un:eq:geod-bounds}. Bedford--Taylor continuity on each cover gives
the homogeneous equation. The comparison principle, proved using
Stokes with the above chart weights, gives uniqueness. It also gives
invariance. Compact families of endpoints have uniform positivity,
subsolutions and the finitely many derivative bounds just listed.
\end{proof}

\begin{lemma}[Convexity and endpoint slopes]\label{un:lem:convex}
The orbifold Mabuchi functional, defined by
\[
d\MM|_u(f)\coloneqq -\int_Z f(S_u-\bar S)\omega_u^n,
\]
is continuous and convex along the geodesics of
Lemma~\ref{un:lem:geod}. At their smooth endpoints it satisfies
\begin{equation}\label{un:eq:slopes}
\MM'(0+)\geq d\MM|_{u_0}(\dot u_0),\qquad
\MM'(1-)\leq d\MM|_{u_1}(\dot u_1).
\end{equation}
If a fixed holomorphic field $Y$ has a real normalized Hamiltonian
$h_u^Y$ on the $T$-invariant subspace $\HH_T$, where
$T\coloneqq \overline{\exp(\R\Impart Y)}$ is compact, then
\[
dE_Y|_u(f)\coloneqq \int_Z f h_u^Y\omega_u^n
\]
defines a functional affine along geodesics in $\HH_T$.
Thus $\MM_Y\coloneqq \MM+E_Y$ has the same convexity and endpoint properties.
\end{lemma}

\begin{proof}
We adapt the local proof of \cite[Section 3]{BB} to the orbifold
charts. The energy-plus-entropy
formula for $\MM$ follows from orbifold Stokes. On a cover write
$\alpha\coloneqq \pi_Z^*\omega+\ddc U=\ddc\Phi$ and let $\Psi$ be the local
logarithm of the horizontal Monge--Ampere density. The convexity
calculation is the positivity of
\[
\ddc\max(\Psi,\chi-a)\wedge\alpha^n,
\]
where $\chi$ is a metric weight on the relative orbifold canonical
bundle with $\ddc\chi\geq-k_0\alpha$. To choose it, let $\chi_Z$
be the local weights of a smooth orbifold Hermitian metric on $K_Z$
and choose $k_0$ so that $\ddc\chi_Z+k_0\omega\geq0$.
Set $\chi_0\coloneqq\pi_Z^*\chi_Z$ and
$\chi\coloneqq\chi_0-k_0U$. Then
\[
\ddc\chi+k_0\alpha
=\pi_Z^*(\ddc\chi_Z+k_0\omega)\geq0.
\]
These are local metric weights; their curvatures agree on overlaps.

Choose nested uniformizing balls with the closure of the smaller
ball contained in the larger one. The potential $\Phi$ is then defined
on a neighborhood of the closed smaller ball. The local weighted
Bergman construction of
\cite[Theorem 2.1, Corollary 2.2 and Theorem 3.3]{BB} applies there.
In particular, Corollary 2.2 supplies the local ball statement, so no
global manifold cover of the orbifold is required.
Its curvature lower bound is
$-k\alpha$; wedging with $\alpha^n$ eliminates this term because
$\alpha^{n+1}=0$. Truncation at $\chi-a$, followed first by the local
Bergman limit and then by $a\to\infty$, gives the asserted positive
current. One may average the approximating positive currents over
$\Gamma$; the limiting current is already invariant. Under a change
of relative canonical frame, both $\Psi$ and $\chi$ acquire the same
pluriharmonic term $\ell$. Hence
\[
\max(\Psi+\ell,\chi+\ell-a)
=\max(\Psi,\chi-a)+\ell.
\]
Applying $\ddc$ therefore gives the same current on overlaps.
The local invariant currents consequently descend and patch on the
orbifold. Orbifold fibre integration, with factor $1/|\Gamma|$,
preserves positivity. This proves distributional convexity.

Pointwise convexity and continuity follow from the localized
truncated-entropy argument of \cite[Theorem 3.4]{BB}. Apply its
regularization on uniformizing balls with an invariant partition of
unity and chart integrals weighted by $1/|\Gamma|$. The bounded
complex Hessian supplies the local estimates; the regularized
localized entropies have a uniform semiconvexity bound and hence
continuous pointwise limits. Canonical-frame changes preserve the
limiting currents, so the pieces define the global functional.
Removing the regularization and truncation preserves convexity.
Convexity and lower semicontinuity of entropy give endpoint continuity.

We next prove the endpoint inequalities directly, following
\cite[Lemma 3.5]{BB}. If
$\rho_t\coloneqq \omega_{u_t}^n/\omega^n$, the supporting-line
inequality for $x\log x$ at the smooth positive $\rho_0$ gives
\[
\operatorname{Ent}(\rho_t)-\operatorname{Ent}(\rho_0)
\geq\int_Z(1+\log\rho_0)(\rho_t-\rho_0)\omega^n.
\]
Divide by $t$ and integrate by parts in the right-hand side.
The $C^1$ convergence at the endpoint and local Bedford--Taylor
continuity justify its limit. Adding the differentiable energy terms
gives the first inequality of \eqref{un:eq:slopes}; reversing time gives
the second.

On $\HH_T$, Stokes gives
$h_u^Y=h_0^Y+Y(u)$ with its integral-zero normalization. The
two-parameter variation calculation of
\cite[Lemma 4.9 and Proposition 4.14]{BB}, performed on the
uniformizing covers, gives
\[
\frac{d^2}{dt^2}E_Y(u_t)
=\int_Z h_{u_t}^Y
  (\ddot u_t-|\bar\partial\dot u_t|_{\omega_{u_t}}^2)
  \omega_{u_t}^n.
\]
The polarized first-variation formula is symmetric, so the defining
one-form is closed on the convex set $\HH_T$. To pass to a weak
geodesic, write $u_{\epsilon,t}\coloneqq U_\epsilon(\cdot,t)$ and set
\[
\alpha_\epsilon\coloneqq\pi_Z^*\omega+\ddc U_\epsilon,
\qquad
Q_\epsilon\coloneqq\ddot u_{\epsilon,t}
-|\bar\partial\dot u_{\epsilon,t}|_{\omega_{u_{\epsilon,t}}}^2.
\]
Since $w=t+\sqrt{-1}s$ and $U_\epsilon$ is independent of $s$,
the block determinant identity and \eqref{un:eq:eps} give
\[
Q_\epsilon\,\omega_{u_{\epsilon,t}}^n
\wedge\sqrt{-1}\,dw\wedge d\bar w
=\frac{4}{n+1}\alpha_\epsilon^{n+1}
=\frac{4\epsilon}{n+1}\underline\Omega^{n+1}.
\] The form $\underline\Omega$ is fixed and smooth, and the Hamiltonians
$h_{u_{\epsilon,t}}^Y$ are uniformly bounded by the $C^1$ estimate.
Fibre integration therefore yields
$|\frac{d^2}{dt^2}E_Y(u_{\epsilon,t})|\leq C\epsilon$.
The $C^1$ convergence and local Bedford--Taylor continuity pass the
functional and its endpoint first variations to the weak limit.
Thus $E_Y$ is affine along the weak geodesic. The unique geodesic
between invariant endpoints remains invariant.
\end{proof}

\begin{lemma}[A strict estimate with a fixed measure]\label{un:lem:strict}
Let $\mu$ be a fixed positive smooth orbifold volume form of mass
$\mathsf V$. Put
\[
I_\mu(u)\coloneqq \int_Zu\mu,\quad
E(u)\coloneqq \frac1{n+1}\sum_{j=0}^n\int_Zu\omega_u^j\wedge\omega^{n-j},
\quad \mathcal F_\mu\coloneqq I_\mu-E.
\]
For a weak geodesic with $\omega_{u_t}\leq C\omega$, there is $c>0$,
depending only on $\omega,\mu,C$, such that
\begin{equation}\label{un:eq:strict}
I_\mu'(1)-I_\mu'(0)\geq
c\left\|u_1-u_0-\frac1{\mathsf V}\int_Z(u_1-u_0)\mu
 \right\|_{L^2(\mu)}^2.
\end{equation}
The functional $E$ is affine on weak geodesics.
\end{lemma}
\begin{proof}
We adapt the fixed-measure estimate of
\cite[Proposition 4.1]{BB}. The positive block matrix of $\alpha$
and its horizontal upper bound
give, almost everywhere, even when that block is degenerate,
\[
\ddot u_t\geq C^{-1}|\bar\partial\dot u_t|_\omega^2.
\]
This is the Schur-complement inequality, or its limit after adding
$\delta\Omega$. Apply the Poincare inequality of the fixed connected
compact orbifold with measure $\mu$:
\[
I_\mu''(t)\geq c\int_Z|\dot u_t-a(t)|^2\mu,
\qquad a(t)\coloneqq \mathsf V^{-1}\int_Z\dot u_t\mu.
\]
The background $W^{2,p}$ regularity from Lemma~\ref{un:lem:geod} justifies
this computation and the endpoint derivatives. Integrate in $t$ and
apply the Hilbert-space Cauchy--Schwarz inequality to obtain
\eqref{un:eq:strict}. Finally $E'=\int\dot u_t\omega_{u_t}^n$ and its
second derivative is the fibre integral of $\alpha^{n+1}$, which is
zero. This also follows by passage through \eqref{un:eq:eps}.
\end{proof}

\subsubsection{Proof of extremal uniqueness}\label{app:sec:uniqueness-proof}

We apply the convex perturbation method of \cite{BB}: minimize a
fixed-measure functional on each automorphism orbit, solve the
linearized equation, and use \eqref{un:eq:strict} to compare the
perturbed metrics. Integrals still use $\omega^n$ without $1/n!$.

\begin{proof}[Proof of Theorem~\ref{un:thm:fixed}]
Use $\omega_0$ as a reference form. Both metrics are invariant under
$T=\overline{\exp(\R\Impart Y)}$, and their potentials can be taken
$T$-invariant. They are critical points of $\MM_Y$ on $\HH_T$, since
$h_u^Y=S_u-\bar S$ at an extremal metric with field $Y$.

Put $G\coloneqq\Aut^0(Z,A,Y)$. The independent calculation in
Lemma~\ref{tc:lem:red} identifies its Lie algebra with the
complexification of the Hamiltonian Killing algebra. Thus $G$ is
reductive and the isometry subgroup $K_i$ of either metric is maximal
compact. Since $q(Z)=0$, connected automorphisms preserve $A$:
its isomorphism class can vary only in $\operatorname{Pic}^0(Z)$.

Set $\mu\coloneqq\omega_0^n$. As in
\cite[Lemma 4.2, Proposition 4.7 and the proof of Theorem 4.16]{BB},
$\mathcal F_\mu$ attains a minimum on each orbit $G\cdot\omega_i$.
Indeed, for $\eta=\omega_i$ and $\nu=\eta^n$, Cartan decomposition
parametrizes the orbit by radial smooth Mabuchi geodesics.
Along them $\mathcal F_\nu'(0)=0$ and
$\mathcal F_\nu''(0)=\int_Z|\bar\partial h|_\eta^2\nu>0$
for nonzero tangent directions modulo constants. Compactness of the
unit sphere and convexity give a uniform linear lower bound along
these rays, hence properness. Writing $\eta=\omega_0+\ddc b$,
orbifold Stokes bounds
$\mathcal F_\mu-\mathcal F_\nu=\int_Zu(\omega_0^n-\eta^n)$
by $\|b\|_\infty$ times fixed positive cohomological integrals.
Thus $\mathcal F_\mu$ is proper as well.
Move each metric by $G$ to a minimum, with potential $\widehat u_i$;
it remains extremal with field $Y$ and $T$-invariant.

At $\widehat u_i$, the differential
$d\mathcal F_\mu(f)=\int f(\mu-\omega_{\widehat u_i}^n)$ annihilates
all real $T$-invariant holomorphy potentials, by the orbit minimum.
It annihilates constants by the equal-mass normalization. On real
$T$-invariant functions the linearization of $d\MM_Y$ at this critical
point is
\[
f\longmapsto L_i f\,\omega_{\widehat u_i}^n,
\qquad L_i\coloneqq \D_i^*\D_i,\qquad
\D_i f\coloneqq \bar\partial\nabla^{1,0}f.
\]
The $Y(f)$ terms in the variations of scalar curvature and $h^Y$
cancel. This operator is real, self-adjoint and elliptic on the
invariant subspace; its kernel is precisely the potentials just
annihilated. Orbifold elliptic theory therefore gives a real smooth
$T$-invariant $v_i$ solving
\begin{equation}\label{un:eq:linear}
L_i v_i=1-\frac{\mu}{\omega_{\widehat u_i}^n}.
\end{equation}
Equivalently one may solve on the orthogonal complement of the kernel
and average over $T$; $L_i$ commutes with this averaging. The invariant
right-hand side pairs with any kernel element just as it pairs with
that element's $T$-average, so it is orthogonal to the full kernel.

For small $s>0$, let $u_i^s\coloneqq \widehat u_i+s v_i$ and
$\MM_{Y,s}\coloneqq \MM_Y+s\mathcal F_\mu$. Equation \eqref{un:eq:linear} gives
\begin{equation}\label{un:eq:residual}
\left|d\MM_{Y,s}|_{u_i^s}(f)\right|
\leq C s^2\|f\|_{C^0}
\end{equation}
for every real invariant $f$. This follows by Taylor expansion of
the smooth density of the first variation.
Join $u_0^s,u_1^s$ by the invariant weak geodesic $u_t^s$.
Its Hessian and endpoint velocities have bounds independent of $s$,
by Lemma~\ref{un:lem:geod}. Convexity and endpoint inequalities give
\begin{align*}
0\leq s\big(I_\mu'(u_t^s)|_{t=1}-I_\mu'(u_t^s)|_{t=0}\big)
&\leq \MM_{Y,s}'(1-)-\MM_{Y,s}'(0+)\\
&\leq d\MM_{Y,s}|_{u_1^s}(\dot u_1^s)
       -d\MM_{Y,s}|_{u_0^s}(\dot u_0^s)\leq C s^2.
\end{align*}
Here we used that $E$ is affine and $\MM_Y$ is convex. Applying
Lemma~\ref{un:lem:strict} yields
\begin{equation}\label{un:eq:closing}
\left\|u_1^s-u_0^s-\mathsf V^{-1}
\int_Z(u_1^s-u_0^s)\mu\right\|_{L^2(\mu)}^2\leq C's.
\end{equation}
Since $u_i^s=\widehat u_i+s v_i$ converges smoothly to $\widehat u_i$
as $s\to0^+$, inequality~\eqref{un:eq:closing} implies that the smooth
function $\widehat u_1-\widehat u_0$ equals its $\mu$-average:
\[
\widehat u_1-\widehat u_0
=\mathsf V^{-1}\int_Z(\widehat u_1-\widehat u_0)\mu
\qquad\text{on }Z.
\]
Consequently $\omega_{\widehat u_1}=\omega_{\widehat u_0}$.
Writing $g_i^*\omega_i=\omega_{\widehat u_i}$ for the two changes by
automorphisms chosen above, the element
$a\coloneqq g_1\circ g_0^{-1}\in G$ therefore satisfies
$a^*\omega_1=\omega_0$, as required.
\end{proof}

\begin{remark}
The orbifold adaptation uses local uniformizing covers, radial
properness and \eqref{un:eq:closing}; no global quotient presentation
is needed. Its constants concern a single fixed orbifold, and its
proof uses no cscK existence or stability theorem.
\end{remark}

We conclude with the finite-group refinement of fixed-orbifold
uniqueness that will be used in Section~\ref{eqv:sec:main}.

\begin{lemma}[Fixed-orbifold uniqueness with finite symmetry]
\label{eqv:lem:fixed}
In Theorem~\ref{un:thm:fixed}, suppose that a finite group $\Gamma$
preserves both metrics and the polarization. Then the isometry can be
chosen in $\operatorname{Cent}_{\Aut^0(Z,A,Y)}(\Gamma)^0$.
\end{lemma}
\begin{proof}
Put $H\coloneqq \Aut^0(Z,A,Y)$, acting effectively on $Z$, and let
$K\coloneqq \operatorname{Stab}_H(\omega_0)$. Theorem~\ref{un:thm:fixed}
places $\omega_1$ in the metric orbit $H/K$, and $K$ is maximal compact.
Cartan decomposition gives a unique $X\in\Lie K$ such that, for the
pushforward action on forms,
\[
\omega_1=(\exp(\sqrt{-1}X))_*\omega_0.
\]
Since $\Gamma$ preserves $\omega_0$ and $Y$, it normalizes $H$ and $K$.
Invariance of $\omega_1$ and uniqueness of these polar coordinates imply
$\operatorname{Ad}_\gamma X=X$ for every $\gamma\in\Gamma$.
Thus $a\coloneqq \exp(\sqrt{-1}X)$ centralizes $\Gamma$, is joined to the identity
in that centralizer, and satisfies $a^*\omega_1=\omega_0$.
\end{proof}

\section{Polarization and Bergman estimates}\label{app:sec:bergman}

We extend the limiting polarization, prove convergence of complete
section spaces, and obtain a Bergman lower bound using a connected
ample divisor away from the exceptional points. We then apply the
same argument to families of polarized cscK surfaces.

\subsection{Polarization and complete section spaces}\label{app:sec:sections}

Following the Hilbert-polynomial argument of Zhang--Zhou
\cite[Section 4]{ZZ}, we combine the numerical assumptions with
polarization convergence to obtain one very ample degree and identify
the geometric and Hilbert limits.

\begin{lemma}[Polarization extraction and extension]\label{tc:lem:line}
The convergence supplied by the compactness property lifts, after a
subsequence and unitary gauges at
$\tau=0$, to convergence of Hermitian holomorphic line bundles on regular
compact subsets. The limit extends over $Z$ as an orbifold line bundle
$L_Z$ with smooth positive Hermitian metric $h_Z$. The gauges may be
chosen independent of $\tau$.
\end{lemma}
\begin{proof}
Choose $\dot\varphi=S-\bar S$ and $h(t)\coloneqq h(0)e^{-\varphi(t)}$, so that
$F_{h(t)}=-\sqrt{-1}\omega(t)$. On a finite good cover of a regular compact
subset, unitary radial gauges express the local connections in terms of
their curvatures. Their derivatives are uniformly bounded, because the
pulled-back curvatures converge smoothly. Transition functions have
bounded derivatives and take values in the compact group $U(1)$.
Subsequences on a nested good-cover exhaustion give a limiting unitary
connection and bundle. Its curvature is $-\sqrt{-1}\omega_Z(0)$.
The limiting $(0,1)$ part is integrable, and so defines a holomorphic
line bundle with its Chern metric.

To extend the bundle, we pull it back from a punctured quotient chart to
$B^4\setminus\{0\}$. It is topologically trivial, since this punctured ball
has $H^2=0$. The limiting curvature extends smoothly to the ball. Take
a smooth invariant primitive for it by radial homotopy. The difference
between this connection and the lifted connection is a closed imaginary
one-form and is exact, because $H^1(B^4\setminus\{0\};\R)=0$.
A unitary gauge therefore gives a smooth connection across the origin.
In this gauge deck transformations preserve the invariant connection;
their fibre factors are constant and define a character
$\Gamma\to U(1)$. This extends the bundle equivariantly. Integrability,
the metric and its positivity also extend. Agreement on punctured
overlaps glues the extensions.

It remains to check that the same identifications work at nearby times.
On the fixed source bundle, the ratio of squared norms is
\begin{equation}\label{tc:eq:metric-ratio}
\frac{h(t_j+\tau)}{h(t_j)}
=\exp\left(-\int_0^\tau(S(t_j+u)-\bar S)\,du\right).
\end{equation}
Its mixed derivatives converge on regular compact subsets. Thus the
same time-independent identifications give convergence for all bounded
$\tau$ intervals. This argument fixes the time-dependent scalar gauge,
not just the curvature.
\end{proof}

\begin{lemma}[No loss of section mass]\label{tc:lem:sections}
For each fixed $k\geq1$ there is a constant $B_k$, independent of time,
such that every $s\in H^0(M,L^k)$ satisfies
\begin{equation}\label{tc:eq:moser}
\sup_M |s|_{h(t)^k}^2\leq B_k\int_M|s|_{h(t)^k}^2dV_t.
\end{equation}
Consequently every ordered orthonormal basis at times $t_j$ has a
subsequence converging smoothly on the regular locus to an orthonormal
set of orbifold holomorphic sections of $L_Z^k$.
\end{lemma}
\begin{proof}
With $\Delta_\omega\coloneqq \operatorname{tr}_\omega\ddc$, the line-bundle
Bochner formula gives $\Delta_\omega|s|^2\geq-2k|s|^2$.
Applying Moser iteration with \eqref{tc:eq:sob} proves
\eqref{tc:eq:moser}. This argument requires no lower bound for the
Ricci curvature of the base. Local elliptic estimates for
the converging Chern connections give smooth subsequential convergence.
The uniform supremum bounds and \eqref{tc:eq:volume} imply that all
$L^2$ inner products are preserved in the limit.
On a local cover of a quotient point, a bounded holomorphic section in
a smooth holomorphic frame extends across the puncture by Hartogs.
It retains its deck equivariance and defines an orbifold section.
\end{proof}

\begin{lemma}[Numerics and no dimension jump]\label{tc:lem:hilbert}
For each limit,
\[
\begin{gathered}
c_1(L_Z)^2=c_1(L)^2,\qquad
K_Z\cdot c_1(L_Z)=K_M\cdot c_1(L),\\
\chi(\mathcal O_Z)=1,\qquad h^1(\mathcal O_Z)=h^2(\mathcal O_Z)=0.
\end{gathered}
\]
In particular $q(Z)\coloneqq h^1(Z,\mathcal O_Z)=0$.
There is an integer $r\geq1$, depending only on the Sobolev bound $C_S$
and the fixed numbers
$d\coloneqq c_1(L)^2$ and $e\coloneqq c_1(M)\cdot c_1(L)$, which works for every
subsequential limit simultaneously. For this $r$, the orbifold bundle
$L_Z^r$ descends to a very ample Cartier line bundle $A_Z$ on the
coarse space, and $A_Z$ and $L^r$ have vanishing higher cohomology
and equal Hilbert polynomials.
At this power the limiting orthonormal set in Lemma~\ref{tc:lem:sections}
is a basis of all of $H^0(Z,A_Z)$.
\end{lemma}
\begin{proof}
We adapt the Hilbert-polynomial argument of Zhang--Zhou
\cite[Lemma 4.3, Theorem 4.4 and Remark 4.7]{ZZ} to Calabi flow
time slices, whose scalar curvatures need not be constant.
Their Theorem 2.1 and the paragraph following it give a uniform
very ample power on the smooth members of their class by Matsusaka's
theorem. Here we also make the power uniform over all possible
orbifold limits, using the normal-variety version of Matsusaka's
theorem after uniformly clearing the local Cartier indices.

\emph{Numerical invariants.}

Volume convergence preserves $c_1(L)^2$. Preservation of $\int S\,dV$ in
Lemma~\ref{tc:lem:extremal}, and the Chern--Weil formula
$\int S\,dV=(2\pi)^2c_1(Z)\cdot c_1(L_Z)$, give the second equality.
Orbifold integrals compute the corresponding rational intersection
numbers on the coarse space.

Each limiting positive orbifold line bundle has an ample Cartier power
on the coarse space, by Baily's orbifold Kodaira--Baily embedding theorem
\cite[Section~7, the ``Global Imbedding Theorem'']{Baily}. Fix this limit temporarily and let $k$ run through sufficiently
large divisible integers. Serre vanishing and Riemann--Roch give
\[
h^0(Z,L_Z^k)
=\frac{k^2}{2}c_1(L)^2-\frac{k}{2}K_M\cdot c_1(L)
+\chi(\mathcal O_Z).
\]
To justify this formula on the singular surface, we resolve $Z$. Quotient
singularities are rational, the pullback of a Cartier divisor has the
same cohomology, and its intersection with exceptional curves is zero.
Thus smooth surface Riemann--Roch gives precisely this formula;
see \cite[equations (4.5)--(4.6)]{ZZ} for this resolution argument.
The injection of Lemma~\ref{tc:lem:sections}, as in
\cite[Corollaries 2.11--2.12]{ZZ}, implies
$\chi(\mathcal O_Z)\geq\chi(\mathcal O_M)=1$.

On the other hand $K_Z\cdot c_1(L_Z)<0$. We use the vanishing
observation of Zhang--Zhou \cite[Remark 4.7]{ZZ}:
a nonzero section of the
dualizing sheaf would give an effective canonical Weil divisor with
negative intersection against an ample divisor, which is impossible.
Serre duality for the Cohen--Macaulay surface gives $h^2(\mathcal O_Z)=0$.
Since $Z$ is connected and normal,
$\chi(\mathcal O_Z)=1-h^1(\mathcal O_Z)\leq1$. Equality follows, and also $h^1(\mathcal O_Z)=h^2(\mathcal O_Z)=0$.
This argument used a Cartier power for each individual limit; it did
not require that power to be uniform.

\emph{A common Cartier power.}
We first obtain a bound for every local group in every limit.
For the volume-ratio lower bound from a homogeneous Sobolev
inequality, see Y.~Wang \cite[Proposition 2.4, p.~1265]{Wang94}.
That proposition is stated for complete noncompact manifolds and
requires no Ricci-curvature assumption. We record the local version
needed for \eqref{tc:eq:sob}, including its lower-order term and
the resulting restriction on the radius.
For a time slice, put $v_{t,x}(s)\coloneqq \Vol_{g(t)}B_{g(t)}(x,s)$.
Apply \eqref{tc:eq:sob} to a distance cutoff equal to one on
$B(x,s/2)$, supported in $B(x,s)$, and with gradient at most $2/s$.
For $0<s\leq V^{1/4}$ this gives
\[
v_{t,x}(s/2)^{1/2}
\leq C_S(4s^{-2}+V^{-1/2})v_{t,x}(s)
\leq5C_Ss^{-2}v_{t,x}(s).
\]
Writing $f_{t,x}(s)\coloneqq s^{-4}v_{t,x}(s)$, we obtain
$f_{t,x}(s)\geq(20C_S)^{-1}f_{t,x}(s/2)^{1/2}$.
Iteration over dyadic radii, followed by the Euclidean small-ball
asymptotics at each fixed smooth time slice, yields
\begin{equation}\label{tc:eq:uniform-noncollapse}
\Vol_{g(t)}B_{g(t)}(x,s)\geq\kappa s^4,
\qquad \kappa\coloneqq (20C_S)^{-2},\qquad 0<s\leq V^{1/4}.
\end{equation}
The constants are independent of $t$ and $x$, and the estimate passes
to every limit by the metric and measure convergence in
Theorem~\ref{app:thm:compactness}. Let $\nu_4$ denote the Euclidean
volume of the unit ball in $\R^4$. At any orbifold point $p\in Z$,
orbifold smoothness and effectiveness give
\[
\lim_{s\downarrow0}s^{-4}\Vol_Z B_{g_Z}(p,s)
=\frac{\nu_4}{|\Gamma_p|}\geq\kappa.
\]
Consequently $|\Gamma_p|\leq Q\coloneqq \lceil\nu_4/\kappa\rceil$ for
every such $p$ and every $Z$. Compare the corresponding group-order
bound in Zhang--Zhou \cite[Theorem 2.3]{ZZ}; the cutoff argument here
uses only our Sobolev bound and compactness property.

Choose once and for all
\begin{equation}\label{tc:eq:common-cartier-power}
q\coloneqq \operatorname{lcm}(1,\ldots,Q).
\end{equation}
The action of $\Gamma_p$ on the fibre of $L_Z$ is a character whose
order divides $|\Gamma_p|$, hence divides $q$. Thus every isotropy
action on $L_Z^q$ is trivial. It follows that $L_Z^q$ descends to an
ordinary line bundle $B_Z$ on the coarse space. Equivalently, for the
coarse-space map $\pi:Z_{\mathrm{orb}}\to Z$,
\[
\pi^*B_Z\simeq L_Z^q.
\]
The ampleness assertion uses positivity in addition to this descent.
The limiting metric $h_Z$ on $L_Z$ has positive orbifold curvature
\[
\sqrt{-1}F_{h_Z}=\omega_Z>0,
\]
so $L_Z^q=\pi^*B_Z$ carries the positive orbifold metric $h_Z^q$ with
curvature $q\omega_Z>0$. Baily's theorem cited above states that a
compact V-manifold carrying a positive V-line bundle admits a projective
embedding by a sufficiently high tensor power of that bundle. Applied to
$L_Z^q$, it gives an integer $m=m(Z)>0$ for which the associated
Kodaira map is an embedding. Since the isotropy action on $L_Z^q$ is
trivial, the sections and the Kodaira map descend to the coarse space;
the descended bundle is $B_Z^m$, so $B_Z^m$ is very ample and $B_Z$ is
ample. The integer $m$ may initially depend on $Z$; uniformity over all
limits is obtained only later from Matsusaka's theorem. Both $q$ and
its divisibility property are independent of the chosen limit.

\emph{A fixed Hilbert polynomial.}
The numerical identities already proved and the resolution argument
above give, for every limit, the same polynomial
\begin{equation}\label{tc:eq:common-hilbert-polynomial}
P_q(m)\coloneqq \chi(Z,B_Z^{\otimes m})
=\frac{q^2d}{2}m^2+\frac{qe}{2}m+1
=\chi(M,L^{qm}).
\end{equation}
This is an equality of Euler characteristics; no uniform vanishing
or Bergman lower bound has been used to establish it.

\emph{A common embedding power.}
Apply Matsusaka's boundedness theorem for normal projective varieties
with rational singularities and a fixed Hilbert polynomial
\cite[Theorem 2.4, pp.~206--207; see also p.~176]{Matsusaka}
to the entire collection of pairs $(Z,B_Z)$, together with $(M,L^q)$.
Its hypotheses hold because quotient singularities are rational,
each $B_Z$ is an ample Cartier line bundle, and the polynomial is
\eqref{tc:eq:common-hilbert-polynomial}. The theorem gives an integer
$m_0=m_0(P_q)\geq1$ such that $B_Z^{m_0}$ and $L^{qm_0}$ are very
ample and, for every $\ell\geq1$ and $i>0$,
\[
H^i(Z,B_Z^{\ell m_0})=0,
\qquad H^i(M,L^{q\ell m_0})=0.
\]
Here $m_0$ includes the preliminary uniform multiple used in the
notation of Matsusaka's theorem. Set
\begin{equation}\label{tc:eq:common-embedding-power}
r\coloneqq q\,m_0(P_q),\qquad A_Z\coloneqq B_Z^{m_0}.
\end{equation}
Thus $A_Z$ is the descent of $L_Z^r$, in agreement with the notation
of Theorem~\ref{app:thm:main}. The Hilbert polynomials of $(Z,A_Z)$
and $(M,L^r)$ are both $P_q(m_0m)$, and vanishing gives
\[
h^0(Z,A_Z)=P_q(m_0)=h^0(M,L^r).
\]
The orthonormal injection from Lemma~\ref{tc:lem:sections} is therefore
surjective. This is the complete-linear-series argument in the first
paragraph of the proof of Zhang--Zhou's Theorem 4.9 \cite{ZZ}.

Finally, $Q$ and $q$ depend only on $C_S$, while $P_q$ depends only
on $q,d,e$; recall also that $V=2\pi^2d$. Matsusaka's integer is chosen
for this fixed polynomial, not for one selected limit. Hence the
single integer $r$ in \eqref{tc:eq:common-embedding-power} can be fixed
before choosing a time sequence and works for every subsequential
limit. Further extraction needed for convergence of bases does not
change $r$. No uniqueness of geometric limits or Bergman lower bound
is used in its construction.
\end{proof}

\begin{lemma}[Identification of Hilbert limits]\label{tc:lem:hilb-limit}
Embed all spaces by the complete bases in Lemma~\ref{tc:lem:hilbert}.
Every Hilbert limit of the embedded $M$ equals the embedded coarse
space $Z$, as a scheme.
\end{lemma}
\begin{proof}
After extraction, write $X_j=\iota_j(M)$ for the embedded time slices,
$W$ for their Hilbert limit, and $B=\iota_\infty(Z)$ for the embedded
reduced coarse space. We retain the three steps in the identification.

\smallskip
\noindent\textit{Step 1: the fixed Hilbert polynomial.}
Lemma~\ref{tc:lem:hilbert} gives
\[
P(k)\coloneqq\chi(M,L^{rk})=\chi(Z,A_Z^k)
=\frac{r^2d}{2}k^2+\frac{re}{2}k+1.
\]
Thus $X_j$, $W$ and $B$ have the same Hilbert polynomial.
The equality for $B$ is numerical and was proved independently
of its identification with $W$.

\smallskip
\noindent\textit{Step 2: scheme-theoretic inclusion.}
Choose $\ell$ sufficiently large for the Grassmannian construction
of the fixed Hilbert scheme and for $\mathcal I_W(\ell)$ to be
globally generated. With
$I_Y(\ell)\coloneqq H^0(\mathbb P^{N-1},\mathcal I_Y(\ell))$,
Hilbert convergence gives $I_{X_j}(\ell)\to I_W(\ell)$.
Every $F\in I_W(\ell)$ is therefore a coefficientwise limit of
$F_j\in I_{X_j}(\ell)$. Complete-basis convergence on the dense
regular locus gives $F|_B=0$. Since $B$ is reduced and
$\mathcal I_W(\ell)$ is globally generated, this implies
$\mathcal I_W\subseteq\mathcal I_B$, hence $B\subseteq W$ as schemes.

\smallskip
\noindent\textit{Step 3: equality from the full polynomial.}
The kernel $\mathcal K$ of $\mathcal O_W\to\mathcal O_B$ has zero
Hilbert polynomial. By Serre vanishing and eventual global generation
\cite[Lemma 30.14.1(1), (4)]{Bertini}, a sufficiently positive twist
of $\mathcal K$ is globally generated and has no sections.
Thus $\mathcal K=0$ and $W=B$.

The full Hilbert polynomial, rather than degree alone, excludes
additional components, multiplicities and embedded structure.
The argument uses regular-locus section convergence, without a
global Bergman lower bound.
\end{proof}

\subsection{Bergman lower bounds}\label{app:sec:partial}

Complete-section convergence gives the Bergman lower bound on regular
compact subsets. Near the exceptional points, we use a limiting ample
divisor avoiding them: connectedness excludes additional zero components,
and the minimum principle gives the bound.

\begin{theorem}[Transfer across finitely many exceptional points]\label{pc:thm:transfer}
Let $(X_j,L_j,h_j,\omega_j)$ be connected smooth compact polarized
K\"ahler surfaces, with
\[
 \sqrt{-1}F_{h_j}=\omega_j.
\]
Let $Z_{\mathrm{orb}}$ be a connected compact K\"ahler orbifold with
isolated quotient singularities, let $Z$ be its normal coarse surface,
and write $\pi:Z_{\mathrm{orb}}\to Z$ for the coarse-space map.
Suppose $L_Z$ is an orbifold holomorphic line bundle with a smooth
Hermitian metric $h_Z$ satisfying
$\sqrt{-1}F_{h_Z}=\omega_Z>0$.

Let $\Sigma\subset Z$ be a finite set containing $\operatorname{Sing}Z$.
Suppose the polarized K\"ahler structures converge smoothly on compact
subsets of $Z\setminus\Sigma$ to
$(L_Z,h_Z,\omega_Z)$. More precisely, the convergence is realized by
smooth embeddings and smooth bundle identifications on neighborhoods
of the retained compact subsets, under which the complex structures,
K\"ahler forms, Hermitian metrics and Chern connections converge in
$C^\infty$. For a retained compact domain with smooth boundary, these
embeddings are defined on a neighborhood of its closure, so they also
identify boundary collars.

Assume there is a fixed integer $k>0$ such that:
\begin{enumerate}
\item $L_Z^k$ descends to a very ample Cartier line bundle $A$ on $Z$,
with $\pi^*A\simeq L_Z^k$;
\item $L^2$-orthonormal bases
$(s_{a,j})_{a=1}^N$ of $H^0(X_j,L_j^k)$ converge in $C^\infty$ on
compact subsets of $Z\setminus\Sigma$, using the same identifications,
to an orthonormal basis $(s_{a,\infty})_{a=1}^N$ of the entire space
$H^0(Z,A)$.
\end{enumerate}
Then $\inf_{X_j}\rho_{k,j}\geq b_k>0$ for all sufficiently large $j$.
No volume convergence is needed for this transfer statement.
\end{theorem}

\begin{lemma}[Connectedness of an ample divisor]\label{pc:lem:connected}
If $X$ is a smooth connected projective surface and $D$ is a nonzero
effective ample Cartier divisor, then $\Supp D$ is connected.
\end{lemma}
\begin{proof}
Otherwise write $D=D_1+D_2$ with nonzero effective divisors of disjoint
support, including their multiplicities. Then $D_1\cdot D_2=0$, whereas
ampleness gives $D_i^2=D\cdot D_i>0$ for both $i$. This contradicts
the Hodge index theorem.
\end{proof}
The lemma applies to the divisor of every nonzero section of $L_j^k$,
since $L_j$ is ample. We do not need very ampleness at this power on
$X_j$, or a generic perturbation of the approximating section.

\begin{lemma}[Persistence of a whole zero component]\label{pc:lem:graph}
Let $s_\infty$ be a section of $A$ whose divisor $D_\infty$ is smooth,
reduced, nonempty and disjoint from $\Sigma$. Suppose $s_j\to s_\infty$
in $C^\infty$ near $D_\infty$. Then, for large $j$, a fixed transferred
tubular neighborhood contains a nonempty compact smooth zero set $D_j^0$
which is open and closed in $\Supp\divisor(s_j)$. It converges smoothly
to $D_\infty$, with multiplicity one.
\end{lemma}
\begin{proof}
The derivative of $s_\infty$ along the real rank-two normal bundle of
$D_\infty$ is an isomorphism onto the underlying real bundle of $A$.
Choose a sufficiently small closed tubular neighborhood $T$ and a smooth
bundle identification with the pullback from $D_\infty$.
The inverse function theorem, uniformly over the compact zero set,
shows that a sufficiently $C^1$-close section has exactly one zero in
each small normal fibre. These zeros form a smooth normal graph.
On the remaining compact part of $T$, including its boundary,
$s_\infty$ is bounded away from zero, so there are no additional zeros
there for large $j$. Higher convergence follows by differentiating the
implicit equation. This identifies the \emph{entire} zero set in $T$.

We now transport $T$ to $X_j$. The graph is compact and is the full zero set in
the interior of this neighborhood; its boundary collar has no zeros.
It is therefore both closed and open in the global zero support.
The section is transverse there, so the divisor multiplicity is one.
The convergence embeddings need only be smooth, not holomorphic: the
graph assertion is a statement about transverse sections of real bundles.
The zero set in $X_j$ is nevertheless a complex divisor because $s_j$
itself is holomorphic.
\end{proof}

\begin{proof}[Proof of Theorem~\ref{pc:thm:transfer}]
Very ampleness, Bertini on $Z\setminus\Sigma$ \cite{Bertini}, and avoidance
of finitely many points give a section $s_\infty$ whose divisor satisfies
Lemma~\ref{pc:lem:graph}. The divisor cannot be empty: an everywhere nonzero
section would trivialize an ample line bundle on a projective surface.
Normalize $\|s_\infty\|_2=1$ and write
\[
s_\infty=\sum_{a=1}^N c_a s_{a,\infty},\qquad
\sum_a|c_a|^2=1,\qquad s_j\coloneqq \sum_{a=1}^N c_a s_{a,j}.
\]
Thus $\|s_j\|_2=1$ and $s_j\to s_\infty$ on all regular compact sets.

By Lemma~\ref{pc:lem:graph}, the zero support contains a nonempty compact
subset $D_j^0$ which is both open and closed. The support is connected
by Lemma~\ref{pc:lem:connected}. Consequently,
\begin{equation}\label{pc:eq:allzeros}
\divisor(s_j)=D_j^0
\end{equation}
as divisors, since the multiplicity along the graph is one. In particular,
the section has no zeros in any region omitted by regular-set convergence.

Choose small mutually disjoint closed orbifold neighborhoods $B_p$ of
the points of $\Sigma$, with smooth boundaries in the regular locus and
disjoint from the tubular neighborhood of $D_\infty$. For ordinary
exceptional points use ordinary balls. Put
\[
K\coloneqq Z\setminus\bigcup_{p\in\Sigma}\operatorname{int}B_p,
\qquad V_j\coloneqq X_j\setminus\operatorname{int}\Phi_j(K).
\]
We use convergence embeddings on a neighborhood of $K$, so that they
also identify its boundary collars. The compact domain $V_j$ may be disconnected,
but every boundary component is in these collars. Every component of
$V_j$ has boundary: a boundaryless component would be open and closed
in connected $X_j$, while $\Phi_j(\operatorname{int}K)$ is nonempty.

Since $s_\infty$ is nonzero on $\partial K$, local convergence gives a
constant $a>0$ such that $|s_j|_{h_j^k}^2\geq a$ on $\partial V_j$.
By \eqref{pc:eq:allzeros}, $s_j$ has no zeros on $V_j$. Hence
\begin{equation}\label{pc:eq:minimum}
u_j\coloneqq \log|s_j|_{h_j^k}^2\text{ is smooth on }V_j,\qquad
\ddc u_j=-k\omega_j,\qquad \Delta_{\omega_j}u_j=-2k<0.
\end{equation}
The minimum principle, applied to each compact component, gives
$u_j\geq\log a$ on $V_j$. Its use needs no uniform boundary geometry or
uniform ellipticity constant across the sequence. By \eqref{pc:eq:bergman},
\begin{equation}\label{pc:eq:badbound}
\rho_{k,j}\geq |s_j|_{h_j^k}^2\geq a\quad\text{on }V_j.
\end{equation}
On $\Phi_j(K)$, the complete-basis convergence gives uniform convergence
of $\rho_{k,j}$ to $\rho_{k,Z}$. Since $A$ is basepoint free,
$\min_K\rho_{k,Z}>0$. Combining this with \eqref{pc:eq:badbound} proves
the result. If $\Sigma$ is empty, convergence on the whole compact limit
already gives the bound and the omitted-region argument is unnecessary.
\end{proof}

\begin{remark}[The use of the minimum principle]
Connectedness first excludes interior zeros, where the logarithm would
be singular. The chosen section controls $V_j$; the complete basis
controls the retained region. The Cartier power is essential: a nontrivial
isotropy character can force every section to vanish at an orbifold point.
\end{remark}

\begin{corollary}[One degree for the entire flow]\label{app:cor:partial}
At the common degree $r$ of Lemma~\ref{tc:lem:hilbert},
$0<b\leq\rho_{r,t}\leq B<\infty$ for all $t\geq0$.
\end{corollary}
\begin{proof}
The upper bound follows from the evaluation characterization
\eqref{pc:eq:bergman} and \eqref{tc:eq:moser}. Fix $r$ as in
\eqref{tc:eq:common-embedding-power} before choosing any time sequence.
If the lower bound failed,
choose $(t_j,x_j)$ with $\rho_{r,t_j}(x_j)\to0$. Since $L^r$ is very ample,
the Bergman function is strictly positive and continuous on every compact
time interval. Thus $t_j\to\infty$ after a subsequence. The compactness property and
Lemmas~\ref{tc:lem:sections}--\ref{tc:lem:hilbert} give the complete-basis
convergence at this already fixed degree. Theorem~\ref{pc:thm:transfer}
gives a constant $b_Z>0$ such that, on the extracted subsequence,
$\inf_M\rho_{r,t_j}\geq b_Z$ for all sufficiently large $j$.
Although $b_Z$ may depend on this limit, its positivity contradicts
$\rho_{r,t_j}(x_j)\to0$. Thus the lower bound for the entire flow is
uniform; equality of the constants obtained from different limits
is not required.
\end{proof}
This proves part (1) of Theorem~\ref{app:thm:main}. The divisor argument
in Theorem~\ref{pc:thm:transfer} does not use the Calabi-flow equation or
the cscK equation. The uniform Sobolev bound and the flow are used in
the preceding compactness and complete-section arguments.
Subsection~\ref{app:sec:alternative} gives an area proof and a uniformization
argument that does not require a common embedding power.

\subsection{Further Bergman estimates and cscK sequences}\label{app:sec:alternative}

We first give an alternative area proof of \eqref{pc:eq:allzeros},
using volume convergence. We then use tensor powers to obtain a
common Bergman degree without a common embedding power, and apply
this to Zhang--Zhou's cscK compactness class
\cite[Problem 8.1]{ZZ}.

\begin{proposition}\label{pc:prop:area}
In the setup of Theorem~\ref{pc:thm:transfer}, suppose also that
$\Vol(X_j,\omega_j)\to\Vol(Z,\omega_Z)$. With the chosen sections and
graphs above, $\divisor(s_j)=D_j^0$ for all large $j$.
\end{proposition}
\begin{proof}
The graph consists of whole multiplicity-one components, so
$E_j\coloneqq\divisor(s_j)-D_j^0$ is an effective integral divisor.
Since $\Vol(X_j)=2\pi^2L_j^2$, volume convergence gives
$L_j^2\to L_Z^2$. Regular convergence along the Cartier divisor
$D_\infty\sim kL_Z$ therefore gives
\[
\int_{E_j}\omega_j
=2\pi kL_j^2-\int_{D_j^0}\omega_j
\longrightarrow2\pi kL_Z^2-\int_{D_\infty}\omega_Z=0.
\]
But integral ampleness on the smooth source implies
$\int_{E_j}\omega_j\in2\pi\Z_{>0}$ whenever $E_j\ne0$.
Thus $E_j=0$ for all large $j$, without assuming integrality
of the orbifold number $L_Z^2$.
\end{proof}

The degree in Theorem~\ref{pc:thm:transfer} may depend on the extracted limit.
The following argument eliminates that dependence without a common
Matsusaka bound.

\begin{lemma}[Tensor powers]\label{pc:lem:tensor}
Suppose, for a fixed $k$, all sections satisfy
$\|s\|_\infty^2\leq B_k\|s\|_2^2$. Then for every integer $m\geq1$,
\begin{equation}\label{pc:eq:tensor}
\rho_{mk}(x)\geq B_k^{1-m}\rho_k(x)^m.
\end{equation}
\end{lemma}
\begin{proof}
If $\rho_k(x)=0$ there is nothing to prove. Otherwise choose a unit section
$s$ attaining evaluation norm at $x$. Then
\[
\|s^m\|_2^2=\int|s|^{2m}dV\leq
\|s\|_\infty^{2m-2}\|s\|_2^2\leq B_k^{m-1}.
\]
Use $s^m$ as a competitor in \eqref{pc:eq:bergman} at degree $mk$.
\end{proof}

\begin{proposition}[Compactness uniformization]\label{pc:prop:uniform}
Let $\mathcal F$ be a family of polarized metrics such that:
\begin{enumerate}
\item for each fixed $k$ the section upper bound of
Lemma~\ref{pc:lem:tensor} holds uniformly in $\mathcal F$;
\item every sequence in $\mathcal F$ has a subsequence satisfying the
assumptions of Theorem~\ref{pc:thm:transfer} for some fixed $k$ along that
subsequence.
\end{enumerate}
Then a single integer $k_0$ and $b>0$ satisfy $\rho_{k_0}\geq b$
on every member of $\mathcal F$.
\end{proposition}
\begin{proof}
If not, the infimum over all members and points is zero at every degree.
For $q_j\coloneqq j!$, choose $(X_j,x_j)\in\mathcal F$ such that
\begin{equation}\label{pc:eq:diagonal}
\rho_{q_j,X_j}(x_j)<\exp(-q_j^2).
\end{equation}
Extract a subsequence and use Theorem~\ref{pc:thm:transfer} to obtain a fixed
$k$ and $b_k>0$ with $\rho_{k,X_j}\geq b_k$. Increase $B_k$ to at least
one. For large original indices $j$, $k$ divides $q_j$; write $q_j=m_jk$.
Lemma~\ref{pc:lem:tensor} yields
\[
\rho_{q_j,X_j}(x_j)\geq
B_k^{1-m_j}b_k^{m_j}\geq\exp(-Cq_j)
\]
with a fixed $C$. This contradicts \eqref{pc:eq:diagonal}.
\end{proof}
If a common very ample degree with complete-basis convergence is already
known for all limits, apply Theorem~\ref{pc:thm:transfer} directly at that
degree and argue by a violating sequence. The factorial argument is
useful because it needs only a degree chosen after extracting the limit.

\begin{corollary}[A uniform Bergman lower bound for cscK surfaces]\label{pc:cor:zz}
Let $\mathcal K(V_0,C_0)$ be the class of polarized cscK surfaces with
volume at most $V_0$, $|S|\leq1$ and scale-invariant Sobolev constant
at most $C_0$. There are an integer $k_0>0$ and a constant $b_0>0$,
depending only on $V_0,C_0$, such that
\[
\inf_X\rho_{k_0,X}\geq b_0^2
\qquad\text{for every }X\in\mathcal K(V_0,C_0),
\]
with the norm convention fixed above.
\end{corollary}
\begin{proof}
Take any sequence $(X_j,L_j,h_j,\omega_j)$ in
$\mathcal K(V_0,C_0)$. By Zhang--Zhou's compactness theorem
\cite[Theorem 2.3]{ZZ}, after passing to a subsequence it converges,
as polarized K\"ahler structures on regular compact subsets, to a
compact polarized cscK orbifold $(Z,L_Z,h_Z,\omega_Z)$. The convergence
embeddings are defined on neighborhoods of those compact subsets.
As in \cite[Section 2.1]{ZZ}, we may also assume that the smooth pairs
$(X_j,L_j)$ have a common Hilbert polynomial.

Choose an integer $q>0$, allowed to depend on this extracted limit,
such that $L_Z^q$ descends to an ample Cartier line bundle $B$ on the
coarse surface $Z$. Zhang--Zhou's Hilbert-polynomial theorem
\cite[Theorem 4.4]{ZZ} gives
\[
 P_q(m)\coloneqq\chi(X_j,L_j^{qm})=\chi(Z,B^m)
\]
for every $m\geq1$ and every $j$ in the extracted sequence.
For the smooth source pairs $(X_j,L_j^q)$, Matsusaka's theorem at this
fixed polynomial supplies an integer $a>0$ such that $L_j^{qa}$ is
very ample and
\[
 H^i(X_j,L_j^{qam})=0
 \qquad(i>0,\ m\geq1)
\]uniformly in $j$; we use the form
\cite[Theorem 2.4]{Matsusaka} already applied in
Lemma~\ref{tc:lem:hilbert}, and compare
\cite[Theorem 2.1 and the following paragraph]{ZZ} for the smooth-source
embedding bound. On the single fixed coarse limit, ampleness of $B$
and Serre vanishing allow us to choose $m>0$ sufficiently large that
\[
 A\coloneqq B^{am}\text{ is very ample},\qquad
 H^i(Z,A)=0\quad(i>0).
\]
Set $k\coloneqq qam$. This degree is fixed along the extracted
subsequence, and the preceding vanishing and polynomial equality give
\[
 h^0(X_j,L_j^k)=P_q(am)=h^0(Z,A).
\]

Zhang--Zhou's section convergence
\cite[Corollary 2.11]{ZZ} now gives, after further extraction,
convergence of orthonormal bases at degree $k$ to an orthonormal set
in $H^0(Z,A)$. Equality of dimensions makes that set a complete basis;
this is the complete-series step in
\cite[the proof of Theorem 4.9]{ZZ}. We have therefore verified every
subsequential hypothesis of Theorem~\ref{pc:thm:transfer} at one
fixed degree along this subsequence. Neither the degree nor the
resulting positive lower bound needs to be chosen before extraction.

Finally, \cite[Proposition 2.10]{ZZ} gives a uniform section upper bound
at each fixed degree, after converting its rescaled norm to our norm.
Since $dV_{k\omega}=k^2dV_\omega$, its bound becomes
$\rho_{k,j}\leq K_0^2k^2$ in our convention, with
$K_0=K_0(V_0,C_0)$.
Proposition~\ref{pc:prop:uniform} applies and yields one degree $k_0$
and one positive lower bound for all of $\mathcal K(V_0,C_0)$.
The class is specified only by $V_0,C_0$, so these constants depend
only on $V_0,C_0$. Renaming the lower bound as $b_0^2$ proves the claim.
\end{proof}

\begin{remark}[Bergman estimates from weak compactness]\label{pc:rem:zz}
Corollary~\ref{pc:cor:zz} answers, in complex dimension two under
the stated bounds, the uniform Bergman lower-bound question in
\cite[Problem 8.1]{ZZ}. It requires neither
\eqref{tc:eq:numerical} nor a Calabi flow. Zhang--Zhou supply
orbifold compactness, complete-section convergence and section
upper bounds; the divisor argument and
Proposition~\ref{pc:prop:uniform} then give one degree and one
positive lower bound for the whole family.

Their rescaled $L^{2,\#}$ norm uses the volume form of $k\omega$,
so $\rho_k^\#=k^{-2}\rho_k$ in dimension two, changing only the
final constant.
\end{remark}

\section{Test configurations and convergence}\label{app:sec:convergence}

We now use the complete section spaces to study the asymptotic dynamics
in a fixed Hilbert scheme. The common limiting generator leads to a
test configuration for each geometric limit. Fixed-orbifold extremal
uniqueness then identifies all limits of the Calabi flow up to polarized
holomorphic isometry. Finally, we use stability to identify this limit
with a smooth cscK metric on the original surface.

\subsection{The common finite-dimensional dynamics}\label{app:sec:dynamics}

We express the asymptotic behaviour of the flow in the
fixed section space obtained in Subsection~\ref{app:sec:sections}. We first
prove convergence of the Gram matrices on bounded time intervals and
identify their limits with the actions of the extremal fields. The
rationality of these actions, together with connectedness, then gives a
common trace-free generator. We use this generator to verify the discrete
hypothesis of Chen--Sun--Wang \cite{CSW} and obtain a connected set of
Hilbert limits in a single orbit closure.

All section spaces now have the fixed dimension
$N=h^0(M,L^r)$. Hermitian pairings are linear in their first entry.
For an ordered column of sections $s$, replacing it by $Ps$ changes its
Gram matrix $G$ to $PGP^*$.
The matrices $A_i$, $C_Z$ and $\Lambda$ below act on the projective
evaluation coordinates $[s_1:\cdots:s_N]$, equivalently on the dual of
the complete section space. The usual induced action on sections is
the dual action and has the opposite weights. Subsection~\ref{app:sec:tc}
uses these usual section weights to define the Donaldson--Futaki invariant.

\begin{lemma}[Convergence of Gram matrices]\label{tc:lem:gram}
Let $s_j$ be any ordered $H_{r,t_j}$-orthonormal bases, where
$H_{r,t}(s,s)=\int_M|s|_{h(t)^r}^2dV_t$.
In fixed source holomorphic sections set
$G_j(\tau)\coloneqq \operatorname{Gram}_{H_{r,t_j+\tau}}(s_j)$.
After subsequence extraction,
\begin{equation}\label{tc:eq:gram}
G_j(\tau)\longrightarrow G_\infty(\tau)=e^{-2\tau C_Z}
\end{equation}
uniformly on bounded intervals, where $C_Z$ is Hermitian. Its projective
action is the holomorphic motion $F_\tau$ in Lemma~\ref{tc:lem:extremal}.
\end{lemma}
\begin{proof}
\noindent\textit{Step 1. Convergence of the global inner products.}
Fix a bounded time interval $I$. After extraction, the ordered bases
$s_j=(s_{j,1},\ldots,s_{j,N})^T$ converge on regular compact subsets
to a complete basis $\boldsymbol\sigma=(\sigma_1,\ldots,\sigma_N)^T$
on $Z$, orthonormal at $\tau=0$. No change of the ordered bases is
needed, so this extraction also applies to the recursion used below.
Let $h_{A,\tau}$ denote the metric on $A_Z$ descended from
$h_Z(\tau)^r$, and set
\[
H_{Z,\tau}(u,v)\coloneqq\int_Zh_{A,\tau}(u,v)\,dV_{Z,\tau}.
\]
The time-independent bundle identifications and
\eqref{tc:eq:metric-ratio} give convergence of the local integrands,
uniformly for $\tau\in I$.

To control their global integrals, choose a regular compact set $K$
whose omitted volumes are less than $(2B_r)^{-1}$, uniformly on $I$
for large $j$, using Lemma~\ref{app:lem:tails}. For any source section
$s$, the section upper bound \eqref{tc:eq:moser} gives
\[
\|s\|_{H_{r,t_j+\tau}}^2
\leq\int_{\Phi_j(K)}|s|_{h(t_j+\tau)^r}^2\,dV_{t_j+\tau}
+\frac12\|s\|_{H_{r,t_j+\tau}}^2.
\]
Absorbing the last term shows that the full squared norm is at most
twice the integral over $\Phi_j(K)$. Local convergence bounds that
integral uniformly for every unit-coefficient linear combination of
the basis sections. In particular, there is $M_I<\infty$ such that
$\|s_{j,a}\|_{H_{r,t_j+\tau}}^2\leq M_I$ for all $a$, all $\tau\in I$
and large $j$. For any larger regular compact set $K'$,
\[
\int_{M\setminus\Phi_j(K')}|s_{j,a}|_{h(t_j+\tau)^r}^2\,dV_{t_j+\tau}
\leq B_rM_I\Vol_{g(t_j+\tau)}(M\setminus\Phi_j(K')).
\]
The right-hand side tends uniformly to zero as $K'$ exhausts the
regular locus, in the sense of \eqref{tc:eq:volume}. Cauchy--Schwarz
gives the same control for mixed pairings. Hence
\[
G_j(\tau)\longrightarrow
G_\infty(\tau)\coloneqq\operatorname{Gram}_{H_{Z,\tau}}(\boldsymbol\sigma)
\quad\text{uniformly on }I,\qquad G_\infty(0)=I.
\]
A diagonal extraction gives this convergence on every bounded interval.

\smallskip
\noindent\textit{Step 2. The pullback group and its infinitesimal generator.}
By Lemma~\ref{tc:lem:extremal}, the limiting flow satisfies
$\omega_Z(\tau)=F_\tau^*\omega_Z(0)$, where $F_\tau$ is a one-parameter
group of holomorphic automorphisms. Lift its generator to $L_Z$ by
the Hamiltonian lift, choosing the vertical constant to agree with
the metric normalization in \eqref{tc:eq:metric-ratio}. The induced
one-parameter group $\widetilde F_\tau$ on $A_Z$ then satisfies
\[
h_{A,\tau}=\widetilde F_\tau^*h_{A,0},\qquad
dV_{Z,\tau}=F_\tau^*dV_{Z,0}.
\]
Define its pullback operator on $H^0(Z,A_Z)$ by
\[
(T_\tau\sigma)(x)\coloneqq
\widetilde F_{\tau,x}^{-1}\bigl(\sigma(F_\tau x)\bigr),
\qquad
\widetilde F_{\tau,x}:(A_Z)_x\longrightarrow(A_Z)_{F_\tau x}.
\]
In the fixed ordered basis, define $P(\tau)$ by
$T_\tau\boldsymbol\sigma=P(\tau)\boldsymbol\sigma$.
The common fibre identification disappears in projective coordinates,
so $P(\tau)$ represents the action of $F_\tau$ on the evaluation
coordinates. Pullback reverses composition, but these automorphisms
form a commuting one-parameter group. Consequently
\[
P(0)=I,\qquad P(\tau+s)=P(\tau)P(s)=P(s)P(\tau).
\]
Define $C_Z\coloneqq\left.\frac{d}{d\tau}\right|_{\tau=0}P(\tau)$.
Differentiating the group identity with respect to $s$ at zero gives
\[
P'(\tau)=P(\tau)C_Z,\qquad P(0)=I,
\]
whose unique solution is $P(\tau)=e^{\tau C_Z}$.
The imaginary-time subgroup acts by Hamiltonian isometries with
a unitary bundle lift. It preserves $h_{A,0}$ and $dV_{Z,0}$, so
its matrices $e^{\sqrt{-1}\theta C_Z}$ are unitary in the orthonormal
basis $\boldsymbol\sigma$. Differentiating
$(e^{\sqrt{-1}\theta C_Z})^*e^{\sqrt{-1}\theta C_Z}=I$ at $\theta=0$
gives $C_Z^*=C_Z$. This uses no periodicity of the imaginary-time
subgroup; rationality is addressed in the next lemma.

\smallskip
\noindent\textit{Step 3. Identification of the limiting Gram matrix.}
The pullback formulas for the metric and volume, followed by change
of variables under $F_\tau$, give
\[
H_{Z,\tau}(T_\tau\sigma_a,T_\tau\sigma_b)
=H_{Z,0}(\sigma_a,\sigma_b)=\delta_{ab}.
\]
Thus the pulled-back basis is orthonormal for the evolving inner
product. Since Hermitian pairings are linear in their first entry,
the same identity in matrix form is
\[
P(\tau)G_\infty(\tau)P(\tau)^*=I.
\]
Solving for $G_\infty(\tau)$ and using $C_Z^*=C_Z$, we obtain
\[
G_\infty(\tau)
=P(\tau)^{-1}\bigl(P(\tau)^{-1}\bigr)^*
=e^{-\tau C_Z}e^{-\tau C_Z}
=e^{-2\tau C_Z}.
\]
Together with Step~1, this proves \eqref{tc:eq:gram}.
\end{proof}

\begin{lemma}[Rational spectrum and a common generator]\label{tc:lem:spectrum}
The trace-free parts of all $C_Z$ have one common spectrum. Thus there
is a fixed trace-free Hermitian matrix $\Lambda$ such that each is
unitarily conjugate to $\Lambda$. A positive multiple of $\Lambda$
has integral eigenvalues.
\end{lemma}
\begin{proof}
We give the coefficient calculation underlying the
Futaki--Mabuchi rationality argument \cite{FM}, including the
contributions of the isolated quotient points.

\smallskip
\noindent\textit{Step 1. Torus weights and Hamiltonian normalization.}
Fix a limit and write $\omega\coloneqq\omega_Z(0)$, with complex
scalar curvature $S$. Put
\[
\omega_A\coloneqq r\omega,\qquad S_A\coloneqq S/r,\qquad
dV_A\coloneqq\omega_A^2/2,\qquad c\coloneqq(2\pi)^{-2}.
\]
Choose a Hamiltonian isometry torus whose Lie algebra contains
$V_{\mathrm{ext}}\coloneqq\frac12J\nabla_\omega S$. Remove
ineffective infinitesimal directions and lift a finite cover of the
torus to $A_Z$; denote the resulting torus by $T$.
To justify this lift, note that $H^1(Z,\mathcal O_Z)=0$ makes the
connected torus preserve $A_Z$. Its inverse image in the linear
stabilizer of the complete section embedding is a central extension
by scalars. A maximal algebraic torus in this inverse image maps onto
the complexified torus; a rational splitting of the cocharacter
lattices gives a lift after a finite cover. On its compact subgroup,
the ratio of the pulled-back metric to the induced metric on $A_Z$
is a positive constant,
since the curvature is preserved. These constants form a character
of a compact group into $\R_{>0}$, so the lifted compact torus is
unitary.
Choose a
cocharacter basis $\lambda_a:\C^*\to T^\C$, $1\leq a\leq\ell$.
Let $\xi_a$ be the real infinitesimal vector field of
$\theta\mapsto\lambda_a(e^{\sqrt{-1}\theta})$, with circle period $2\pi$.
We use the usual section action
$\rho_k(g)s\coloneqq g\cdot s(g^{-1}\cdot)$ and define its
commuting Hermitian weight operators by
\[
\rho_k(\lambda_a(e^{\sqrt{-1}\theta}))
       =\exp(\sqrt{-1}\theta D_{a,k}).
\]
Thus the simultaneous eigenvalue vectors of the $D_{a,k}$ belong
to $\Z^\ell$.

Let $f_a$ be the Hamiltonian with its additive constant fixed by
this bundle lift. With the usual section action the convention is
\[
\iota_{\xi_a}\omega_A=df_a.
\]
At a torus-fixed point, $f_a$ equals the weight of the action on
the fibre of $A_Z$. These Hamiltonians need not have zero mean;
we will center them after computing the coefficients.

\smallskip
\noindent\textit{Step 2. The equivariant orbifold Riemann--Roch formula.}
Set $\alpha\coloneqq\omega_A/(2\pi)$ and
$\beta\coloneqq\rho_A/(2\pi)$, where $\rho_A$ is the Ricci form.
For $\mathbf t=(t_1,\ldots,t_\ell)$ near zero, put
$g(\mathbf t)\coloneqq\prod_a\lambda_a(e^{t_a})$ and define the
equivariant Euler characteristic using the induced actions on cohomology:
\[
\chi_k(\mathbf t)\coloneqq
\sum_{j=0}^2(-1)^j\tr\!\left(
 \rho_k(g(\mathbf t))
 \bigm|H^j(Z,A_Z^k)\right).
\]
For $k$ sufficiently large, higher-cohomology vanishing identifies
this with $\tr\exp(\sum_a t_aD_{a,k})$.
In variables normalized so that a weight $\mu$ contributes
$e^{\langle\mu,\mathbf t\rangle}$, the equivariant first Chern
forms are
\[
c_1^T(A_Z)(\mathbf t)=\alpha+\sum_a t_af_a,\qquad
c_1^T(TZ)(\mathbf t)=\beta+\sum_a t_a\tau_a,\qquad
\tau_a=-\Delta_A f_a.
\]
Here $TZ$ denotes the orbifold tangent bundle and
$\Delta_A f\coloneqq g_A^{i\bar j}\partial_i\partial_{\bar j}f$.
The formula for $\tau_a$ comes from the natural lift to
$\det TZ$: the complex divergence of
$\xi_a^{1,0}=-\sqrt{-1}\,g_A^{i\bar j}(f_a)_{\bar j}\partial_i$
is $-\sqrt{-1}\Delta_A f_a$. In particular,
$\int_Z\tau_a\,dV_A=0$ by orbifold integration by parts.
If $\Theta_T(\mathbf t)$ is the normalized equivariant curvature
of the Chern connection on $TZ$, the equivariant Todd form means
\[
\operatorname{Td}_T(TZ)(\mathbf t)
\coloneqq\det\!\left(\frac{\Theta_T(\mathbf t)}
                          {I-e^{-\Theta_T(\mathbf t)}}\right),
\qquad
\tr\Theta_T(\mathbf t)=\beta+\sum_a t_a\tau_a,
\]
with the quotient interpreted by its power series at zero.

For a singular point $p$, let $\Gamma_p$ be its local group.
The connected torus fixes $p$. We choose the finite cover in Step~1
also so that its local lifts commute with every $\Gamma_p$.
The group of local lifts at a quotient point is a finite extension of
the torus. Its identity component covers the torus after an isogeny
and commutes with the finite chart group, because its conjugation
action on that group is continuous and hence constant. Taking a
common finite cover for the finitely many quotient points gives the
stated choice.
Let
$E_p(\mathbf t)$ be its complexified action on the local tangent
space and $\ell_p$ its weight vector on $(A_Z)_p$. Thus
$E_p(0)=I$ and $\ell_{p,a}=f_a(p)$.
The formula we use is
\begin{equation}\label{tc:eq:equivariant-rr}
\begin{aligned}
\chi_k(\mathbf t)
={}&\int_Z e^{k(\alpha+\sum_a t_af_a)}
                    \operatorname{Td}_T(TZ)(\mathbf t)\\
&+\sum_{p\in\operatorname{Sing}Z}
  \ \sum_{\substack{[g]\in\operatorname{Conj}(\Gamma_p)\\g\ne1}}
 \frac{e^{k\langle\ell_p,\mathbf t\rangle}}
 {|C_{\Gamma_p}(g)|\,
       \det\!\left(I-g^{-1}E_p(\mathbf t)^{-1}\right)} .
\end{aligned}
\end{equation}
The first integral is an orbifold integral and selects ordinary
differential-form degree four. The second sum runs over nonidentity
conjugacy classes; $C_{\Gamma_p}(g)$ is the centralizer of $g$.
There is no additional isotropy character in the numerator:
local isotropy acts trivially on $A_Z$, since it is Cartier.
This is Vergne's equivariant index theorem
\cite[Theorem 1, p.~643]{Vergne}, in the Riemann--Roch form of
\cite{Meinrenken}, specialized to isolated quotient singularities.
Setting $\mathbf t=0$ gives Kawasaki's
nonequivariant formula
\cite[unnumbered theorem, formula (IV), p.~153]{Kawasaki}.

For $g\ne1$, the local action is free away from the origin, so
$g$ has no tangent eigenvalue $1$. Each denominator in the second
line of \eqref{tc:eq:equivariant-rr} is therefore nonzero near
$\mathbf t=0$. Its reciprocal is analytic and independent of $k$.
Consequently, an order-$j$ derivative of a nonidentity term at
zero is a polynomial in $k$ of degree at most $j$. For $j=0,1,2$
these terms are respectively $O(1)$, $O(k)$ and $O(k^2)$.
This also applies to mixed second derivatives.

\smallskip
\noindent\textit{Step 3. Comparison of the powers of $k$.}
Differentiating the character at zero gives
\[
\begin{aligned}
d_k&\coloneqq\chi_k(0)=\dim H^0(Z,A_Z^k),\\
w_{a,k}&\coloneqq\partial_{t_a}\chi_k(0)=\tr D_{a,k},\\
Q_{ab,k}&\coloneqq\partial_{t_a}\partial_{t_b}\chi_k(0)
                         =\tr(D_{a,k}D_{b,k}).
\end{aligned}
\]
Write their leading terms as
\[
d_k=a_0k^2+a_1k+O(1),\qquad
w_{a,k}=b_{0a}k^3+b_{1a}k^2+O(k),\qquad
Q_{ab,k}=d_{ab}k^4+O(k^3).
\]
To determine these coefficients, use
\[
e^{k\alpha}=1+k\alpha+\frac{k^2\alpha^2}{2},
\qquad
\operatorname{Td}_T(TZ)
=1+\frac12\left(\beta+\sum_a t_a\tau_a\right)
       +\text{terms of equivariant degree at least four}.
\]
The first equality is understood after discarding forms of degree
greater than four. Each $t_a$ has equivariant degree two.

At $\mathbf t=0$, the $k^2$ and $k$ terms of ordinary form
degree four are respectively $k^2\alpha^2/2$ and
$k\alpha\beta/2$. Therefore
\[
a_0=\int_Z\frac{\alpha^2}{2},\qquad
a_1=\frac12\int_Z\alpha\beta.
\]
For one derivative, differentiating the exponential contributes
$kf_a$. Its products with $k^2\alpha^2/2$ and
$k\alpha\beta/2$ give the $k^3$ term and one $k^2$ term.
Differentiating the Todd form contributes another $k^2$ term,
namely $(\tau_a/2)(k^2\alpha^2/2)$. Hence
\[
b_{0a}=\int_Zf_a\frac{\alpha^2}{2},\qquad
b_{1a}=\frac12\int_Zf_a\alpha\beta
              +\frac12\int_Z\tau_a\frac{\alpha^2}{2}.
\]
The last integral is zero, since
$\alpha^2/2=c\,dV_A$ and $\int_Z\Delta_A f_a\,dV_A=0$.
The derivatives of the higher Todd terms give at most $O(k)$
here, by ordinary form degree.

For two derivatives, the $k^4$ term comes only from
$k^2f_af_b$ times $k^2\alpha^2/2$. If either derivative falls
on the Todd form, the resulting term has degree at most three
in $k$. Thus
\[
d_{ab}=\int_Zf_af_b\frac{\alpha^2}{2}.
\]
The nonidentity-sector orders in Step~2 do not affect any of
these coefficients. Finally,
$\alpha\beta=c\,S_A\,dV_A$, because $S_A$ is complex scalar
curvature. The coefficient identities are therefore
\[
\begin{aligned}
a_0&=c\int_ZdV_A,&
a_1&=\frac c2\int_ZS_A\,dV_A,\\
b_{0a}&=c\int_Zf_a\,dV_A,&
b_{1a}&=\frac c2\int_Zf_aS_A\,dV_A,\\
d_{ab}&=c\int_Zf_af_b\,dV_A.&&
\end{aligned}
\]

These are coefficients of exact polynomials for all sufficiently
large integers $k$, not merely coefficients in asymptotic expansions.
Indeed, after at most two torus derivatives, the identity-sector
integral in \eqref{tc:eq:equivariant-rr} is a polynomial in $k$
by form-degree truncation; Step~2 gives the same assertion for
every remaining term. The integers $d_k$, $w_{a,k}$ and $Q_{ab,k}$
are respectively dimensions, sums of integral weights, and sums of
products of simultaneous integral weights. A polynomial taking
integer values at all sufficiently large integers has rational
coefficients: evaluate it at sufficiently many consecutive integers
and invert the resulting integer Vandermonde matrix. Consequently
all the displayed coefficients $a_0,a_1,b_{0a},b_{1a},d_{ab}$
are rational.

\smallskip
\noindent\textit{Step 4. The rational linear system for the extremal field.}
The preceding identities give
\[
\bar f_a\coloneqq\frac{\int_Zf_a\,dV_A}{\int_ZdV_A}
                  =\frac{b_{0a}}{a_0},\qquad
\bar S_A=\frac{2a_1}{a_0},\qquad
u_a\coloneqq f_a-\bar f_a.
\]
Define the centered pairing and Futaki vector by
\[
\begin{aligned}
B_{ab}&\coloneqq c\int_Zu_au_b\,dV_A
             =d_{ab}-\frac{b_{0a}b_{0b}}{a_0},\\
F_a&\coloneqq c\int_Zu_a(S_A-\bar S_A)\,dV_A
             =2\left(b_{1a}-\frac{a_1b_{0a}}{a_0}\right).
\end{aligned}
\]
Thus $B$ and $F$ are rational. The matrix $B$ is positive
definite: a zero norm would make a linear combination of the $u_a$
identically zero, and its Hamiltonian vector field would be an
ineffective infinitesimal direction.

Write $V_{\mathrm{ext}}=\sum_a v_a\xi_a$. Since
$\iota_{V_{\mathrm{ext}}}\omega_A=-\frac r2\,dS$, taking
mean-zero Hamiltonians gives
\[
\sum_a v_au_a=-\frac{r^2}{2}(S_A-\bar S_A),\qquad
Bv=-\frac{r^2}{2}F.
\]
Hence $v=-\frac{r^2}{2}B^{-1}F\in\Q^\ell$.
If the torus is trivial, $V_{\mathrm{ext}}=0$ and the same
rationality conclusion is immediate.

The operator in Lemma~\ref{tc:lem:gram} is the pullback action,
which is inverse to the usual section action used above. Any
difference between the two chosen lifts of $V_{\mathrm{ext}}$
is scalar and disappears on taking trace-free parts. In a
simultaneous orthonormal weight basis this gives
\[
(C_Z)_0=-\sum_a v_a(D_{a,1})_0.
\]
The $D_{a,1}$ commute and have simultaneous integral eigenvalues;
their trace-free eigenvalues are rational. It follows that every
eigenvalue of $(C_Z)_0$ is rational in the conventions of
\eqref{app:eq:conventions}. In particular, all limits give ordered
spectra in the same countable set $\Q^N$.

\smallskip
\noindent\textit{Step 5. Independence of the subsequence.}
For real $t\geq0$, the ordered spectrum
\[
\sigma(t)\coloneqq \operatorname{ordered\,spec}
\left(-\tfrac12\log\operatorname{Gram}_{H_{r,t+1}}(s_t)\right)_0
\]
is independent of the $H_{r,t}$-orthonormal basis $s_t$ and is continuous.
The subscript $0$ denotes removal of the scalar part.
Lemma~\ref{tc:lem:gram} makes its tail relatively compact, with every
cluster value one of the rational spectra just described. The cluster
set of a relatively compact continuous curve on a half-line is connected.
A nonempty connected subset of a countable subset of $\R^N$ consists
of a single point.
This proves the common-spectrum assertion. Choose $\Lambda$ to be
the diagonal matrix with this ordered spectrum. Clearing its
finitely many rational denominators gives a positive integer $b$
such that $b\Lambda$ has integral eigenvalues.
\end{proof}

We now choose compatible bases at consecutive integer times. Starting
from an orthonormal basis $s_0$ at time zero, define recursively
\begin{equation}\label{tc:eq:recursion}
Q_i\coloneqq \operatorname{Gram}_{H_{r,i+1}}(s_i),\qquad
s_{i+1}\coloneqq Q_i^{-1/2}s_i,\qquad s_i=A_i s_0.
\end{equation}
Then $B_{i+1}\coloneqq A_{i+1}A_i^{-1}=Q_i^{-1/2}$. Put
\[
\widehat A_i\coloneqq (\det A_i)^{-1/N}A_i,\qquad
\widehat B_{i+1}\coloneqq \widehat A_{i+1}\widehat A_i^{-1}
               =(\det B_{i+1})^{-1/N}B_{i+1}.
\]
The determinants are positive by the recursion. The bases $s_i$ remain
orthonormal and unchanged; this normalization of the coordinate matrices
does not change any projective embedding.
For every subsequence $i_j$, Lemma~\ref{tc:lem:gram} gives in the same basis
\[
G_j(1)\to e^{-2C},\qquad G_j(2)\to e^{-4C}.
\]
The Gram matrix for the second step in its updated basis is
$G_j(1)^{-1/2}G_j(2)G_j(1)^{-1/2}\to e^{-2C}$.
Consequently, after determinant normalization,
\begin{equation}\label{tc:eq:star}
\widehat B_{i_j+1}\to ue^\Lambda u^{-1},\qquad
\widehat B_{i_j+2}\to ue^\Lambda u^{-1}
\quad\text{for some }u\in U(N).
\end{equation}
Thus the matrices $\widehat A_i$ satisfy the discrete hypothesis of
\cite[Section 2]{CSW}.
The verification uses only uniform convergence of Gram matrices and
does not require convergence of their derivatives or of a continuous
horizontal lift.

\begin{lemma}[The compatible Hilbert limit set]\label{app:lem:limitset}
Let $G_\Lambda$ be the centralizer of $\Lambda$ in $GL(N,\C)$ and
$K_\Lambda\coloneqq G_\Lambda\cap U(N)$. There is a nonempty compact connected
$K_\Lambda$-invariant Hilbert limit set $\mathscr C$, every point of which
is represented by a geometric orbifold limit with a complete orthonormal
basis. Each $[Z]\in\mathscr C$ is fixed by $\Lambda$, which induces its
extremal flow. There are $C_0\in GL(N,\C)$ and
$[\overline M]\coloneqq \lim_{t\to\infty}e^{t\Lambda}C_0[M]$ such that
\begin{equation}\label{un:eq:common-orbit}
\mathscr C\subset\overline{G_\Lambda\cdot[\overline M]}.
\end{equation}
Every integer-time geometric limit is represented in $\mathscr C$ up to
a unitary change of the complete section basis.
\end{lemma}
\begin{proof}
Place the fixed Hilbert scheme in a finite-dimensional equivariant
projective representation. Condition \eqref{tc:eq:star} is CSW's discrete
condition $(*)$. Choose unitary gauges $u_i$ as in its equivalent
formulation $(*)'$ and let $\mathscr C$ be the union of the
$K_\Lambda$-orbits of all subsequential limits of
$u_i^{-1}\widehat A_i[M]$. This set is compact and connected by
\cite[Lemma 2.6, pp.~3150--3151]{CSW}.

We describe the construction of $C_0$. On the evaluation-coordinate
space $E\coloneqq\C^N$, let $\lambda_1>\cdots>\lambda_\ell$ be the
distinct eigenvalues of $\Lambda$, with eigenspaces $U_j$.
The growth rates and the resulting filtration are
\[
d(v)\coloneqq\lim_{i\to\infty}\frac1i\log\|\widehat A_i v\|,
\qquad
F_j\coloneqq\{0\}\cup\{v\ne0:d(v)\leq\lambda_j\},
\]
\[
E=F_1\supset\cdots\supset F_\ell\supset F_{\ell+1}=\{0\}.
\]
Existence of $d(v)$ is \cite[Lemma 2.3, p.~3150]{CSW}.
Proposition 2.7 of that paper, on p.~3151, supplies subspaces $W_j$
such that
\[
F_j=W_j\oplus F_{j+1},\qquad
u_i^{-1}\widehat A_iW_j\longrightarrow U_j
\quad\text{in the Grassmannian}.
\]
In particular, $\dim W_j=\dim U_j$. Choose linear isomorphisms
$W_j\to U_j$ and combine them to obtain $C_0\in GL(E)$ with
\[
C_0(W_j)=U_j,\qquad
C_0(F_j)=\bigoplus_{k\geq j}U_k.
\]
Thus $C_0$ identifies the filtration determined by the flow with the
weight filtration of the fixed matrix $\Lambda$. The eigenvalues
alone do not identify these subspaces in the initial coordinates.

The near-identity corrections $R_i\to I$ and equation (7) in
\cite[p.~3153]{CSW} place the gauged sequence in $G_\Lambda$
without changing its limits. Proposition 2.10 there gives
\[
[\overline M]\coloneqq\lim_{t\to\infty}e^{t\Lambda}C_0[M],
\qquad
\mathscr C\subset\overline{G_\Lambda\cdot[\overline M]},
\]
which is \eqref{un:eq:common-orbit}.

The matrix $C_0$ is chosen once from the filtration of the full
sequence, independently of the subsequence defining a geometric
limit. Its existence imposes no additional hypothesis. One may
absorb it into the notation for the initial embedding by replacing
$[M]$ with $C_0[M]$, which represents the same polarized surface.
However, this changes the already fixed Hilbert point, and the
corresponding section basis need not remain orthonormal. We therefore
keep $C_0$ explicit; it cannot in general be replaced by the identity
while both $[M]$ and $\Lambda$ remain fixed. Although the splitting
and $C_0$ involve choices, the intermediate orbit
$G_\Lambda\cdot[\overline M]$ is independent of them
\cite[Proposition 2.11, p.~3154]{CSW}.

Finally, Lemma~\ref{tc:lem:hilb-limit} identifies the gauged Hilbert
limits with the reduced coarse spaces of geometric orbifold limits.
The actual section bases undergo only the unitary gauges $u_i$ and
remain orthonormal; the $R_i$ above serve solely for the algebraic
comparison of limits. Lemma~\ref{tc:lem:gram} identifies the action
of $\Lambda$ with the limiting extremal flow. Every original
integer-time limit is recovered up to a unitary basis change by
passing to a convergent subsequence of the $u_i$. These conclusions
do not use uniqueness of the geometric limit.
\end{proof}

\subsection{A test configuration for each limit}\label{app:sec:tc}

We construct a test configuration whose central fibre is the polarized
variety underlying any prescribed geometric limit. We use the orbit closure
obtained in Subsection~\ref{app:sec:dynamics}, together with reductivity of
the extremal centralizer, to reach the limit by commuting algebraic
one-parameter subgroups. A weight comparison combines these subgroups
into one test configuration. We then compute its Donaldson--Futaki
invariant and show that it is negative for a non-cscK limit and zero
for a cscK limit. The construction identifies the central polarized
variety; it need not induce the flow's weighted filtration, and no
optimality is asserted. We explain this distinction, and its relation to the
K\"ahler--Ricci-flow construction, in
Remark~\ref{tc:rem:two-stage} below.

\begin{lemma}[Reductivity of the extremal centralizer]\label{tc:lem:red}
Let $(Z,A,\omega)$ be a polarized compact orbifold as above with an
extremal metric, and let $Y$ be its extremal holomorphic field.
The centralizer $\Aut(Z,A,Y)$ in the polarized automorphism group
$\Aut(Z,A)$, defined in Subsection~\ref{app:sec:orbifold-uniqueness},
is reductive. In particular $\Aut(Z,A)$ is reductive when $\omega$
is cscK.
\end{lemma}
\begin{proof}
We work in the complete section embedding. Polarized holomorphic vector fields
lift to $A$; their Chern potentials $f$ satisfy
$X=\nabla^{1,0}f$, with a fixed harmless convention factor. They are smooth
orbifold functions. Indeed, vector fields on the coarse quotient lift
on the punctured finite cover and extend by Hartogs.

Put $K\coloneqq J\nabla_{\mathrm{real}}S$. A holomorphic field in the indicated
centralizer commutes with $K$. Since $K$ preserves $\omega$, its
mean-zero Chern potential is $K$-invariant. For
$\mathcal Df\coloneqq \bar\partial\nabla^{1,0}f$, the K\"ahler identity for
$\mathcal D^*\mathcal D$ says that its difference from its conjugate is
a fixed imaginary multiple of differentiation along $K$.
Thus on $K$-invariant functions this operator is real. If
$\mathcal Df=0$, its real and imaginary parts separately lie in the
kernel of $\mathcal D^*\mathcal D$, and hence in the kernel of
$\mathcal D$. For a real kernel function the Hamiltonian vector field
is Killing and real holomorphic. It commutes with $K$, and therefore
with the extremal holomorphic field.

All adjoint identities follow by integration on finite smooth charts
and an invariant partition of unity; the compact orbifold has no
boundary term. We have proved that the complex Lie algebra of the
centralizer is the complexification of its compact Hamiltonian isometry
algebra. In the complete section representation these compact
transformations are unitary. The corresponding linear algebraic
centralizer therefore has reductive identity component, and is reductive
(its component group is finite). The same argument with $K=0$ proves
the cscK assertion.
\end{proof}

Lemma~\ref{app:lem:limitset} gives the common intermediate degeneration
and the required orbit inclusion. We also use the accessibility result
\cite[Proposition 1(2), Section 5.2]{Donaldson}: a boundary point whose
stabilizer is reductive is reached by an algebraic one-parameter subgroup
from a suitable point of the given orbit.
The slice argument applies to the reductive centralizer $G_\Lambda$:
it uses reductivity of the stabilizer of the limiting point and the
affine Hilbert--Mumford criterion in the slice.
Donaldson's slice proof in Section 5.2 applies without assuming that
the source stabilizer is finite.

\begin{lemma}[Combining commuting degenerations]\label{tc:lem:combine}
Let $\lambda,\mu:\C^*\to\mathrm{GL}(V)$ be commuting algebraic
one-parameter subgroups, where $V$ is a finite-dimensional complex
vector space. Suppose $x,y,p\in\mathbf P(V)$ satisfy
$\lim_{z\to0}\lambda(z)x=y$ and $\lim_{z\to0}\mu(z)y=p$.
For every sufficiently large positive integer $m$, the algebraic
one-parameter subgroup $\nu_m(z)\coloneqq\lambda(z)^m\mu(z)$ satisfies
\[
\lim_{z\to0}\nu_m(z)x=p.
\]
\end{lemma}
\begin{proof}
Simultaneously decompose $V=\bigoplus_{(a,b)}V_{a,b}$ into the
integral weight spaces of $\lambda$ and $\mu$, and write a
representative of $x$ as $v=\sum v_{a,b}$. Let $(a_0,b_0)$ be the
lexicographically smallest pair with $v_{a,b}\ne0$. The two
successive projective limits give $p=[v_{a_0,b_0}]$.
Since only finitely many weights occur, every sufficiently large
integer $m$ satisfies
\[
m(a-a_0)+(b-b_0)>0
\quad\text{if }v_{a,b}\ne0,\ (a,b)\ne(a_0,b_0).
\]
The action of $\nu_m$ on $V_{a,b}$ has weight $ma+b$, so
$z^{-(ma_0+b_0)}\nu_m(z)v\to v_{a_0,b_0}$.
This proves equality with the specified projective point.
\end{proof}

\begin{proposition}[An ordinary configuration to each limit]\label{tc:prop:tc}
Every polarized geometric limit $(Z,A_Z)$ is the central fibre of a
single ordinary test configuration for $(M,L)$ of exponent $r$.
The total space is normal.
\end{proposition}
\begin{proof}
Once the preceding lemmas provide an algebraic one-parameter subgroup
whose Hilbert limit is the polarized geometric limit, the construction
of the test configuration is standard. We include the details to
identify its central fibre and polarization precisely and to verify
normality of the total space.

\textit{Step 1: construct a single algebraic degeneration.}
We first consider a limit at integer times, embedded in coordinates
compatible with \eqref{tc:eq:star}. Choose $b>0$ such that $b\Lambda$
has integral eigenvalues, and put
\begin{equation}\label{tc:eq:lambda}
\lambda(z)\coloneqq \exp(-b\Lambda\log z).
\end{equation}
In an eigenbasis of $\Lambda$, the entries of $\lambda(z)$ are integral
powers of $z$, so this is a single-valued algebraic one-parameter subgroup.
The negative sign ensures that $z\to0$ corresponds to forward time
in the asymptotic model of the Calabi flow.
Then $\lambda$ degenerates $C_0[M]$ to $[\bar M]$ and
fixes $[Z]$.

The stabilizer of $[Z]$ in $G_\Lambda$ is reductive by
Lemma~\ref{tc:lem:red}; taking the full linear stabilizer only adds central
scalars or a finite cover. Accessibility gives $g\in G_\Lambda$ and
a one-parameter subgroup $\mu$ of $G_\Lambda$ with
$\lim_{z\to0}\mu(z)g[\bar M]=[Z]$.
If $[Z]$ already belongs to $G_\Lambda[\bar M]$, choose $g$ with
$g[\bar M]=[Z]$ and take $\mu$ to be trivial.
Because $\lambda$ is central in $G_\Lambda$, it commutes with both $g$
and $\mu$. Apply Lemma~\ref{tc:lem:combine} with
$x=gC_0[M]$, $y=g[\bar M]$, and $p=[Z]$ in the fixed equivariant
projective embedding of the Hilbert scheme. It gives the equality
of Hilbert points
\begin{equation}\label{tc:eq:nu}
\nu_m(z)\coloneqq \lambda(z)^m\mu(z),\qquad
\lim_{z\to0}\nu_m(z)gC_0[M]=[Z]
\end{equation}
for all sufficiently large integers $m$.
If $\Lambda=0$, the first degeneration is trivial and accessibility
applies directly to $[Z]\in\overline{G_\Lambda[M]}$.

\textit{Step 2: the polarized family and its normality.}
Put $M'\coloneqq gC_0\iota(M)$, where $\iota$ is the complete
$L^r$-embedding. For the common Hilbert polynomial $P$, let
$\mathcal H\coloneqq\operatorname{Hilb}_P(\mathbb P^{N-1})$.
The orbit map $z\mapsto\nu_m(z)[M']$ extends to
$f:\mathbb A^1\to\mathcal H$ by projectivity, with $f(0)=[Z]$.
Pulling back the universal Hilbert family $\mathcal U$ gives
\[
\mathcal X\coloneqq f^*\mathcal U
\subset\mathbb P^{N-1}\times\mathbb A^1,\qquad
\mathcal L\coloneqq
\operatorname{pr}_1^*\mathcal O_{\mathbb P^{N-1}}(1)|_{\mathcal X}.
\]
This is flat and projective, with polarized fibres $(M,L^r)$ away
from zero and exactly $(Z,A_Z)$ at zero, including the scheme
structure identified in Lemma~\ref{tc:lem:hilb-limit}; compare the
Hilbert-arc construction in \cite[Section 2.2]{Donaldson}.

Flatness makes the restriction over $\C^*$ schematically dense,
so $\mathcal X$ is the scheme-theoretic orbit closure. Consequently
\[
a\cdot(q,z)\coloneqq(\nu_m(a)q,az)
\]
preserves $\mathcal X$, and the natural linearization of
$\mathcal O(1)$ induces one on the relatively very ample bundle
$\mathcal L$. Over $\C^*$ the equivariant trivialization is
$(p,z)\mapsto(\nu_m(z)gC_0\iota(p),z)$, with trivial action on $M$;
it identifies the line bundle with $\operatorname{pr}_M^*L^r$.
Thus this is a test configuration of exponent $r$, independently
of the integer $m$.

Every fibre is normal: the noncentral fibres, including the generic
fibre, are smooth, and the central fibre is the normal surface $Z$.
The base $\mathbb A^1$ is normal, so flatness gives normality of
$\mathcal X$ by \cite[Lemma 10.163.8, Tag 0C22]{Bertini}.
No normalization or change of central fibre is required.

\textit{Step 3: include limits at arbitrary times.}
Write an arbitrary sequence as $t_j=n_j+\tau_j$, with
$\tau_j\in[0,1)$ and, after a subsequence, $\tau_j\to\tau_*$. The compactness property at
$n_j$ and Lemma~\ref{tc:lem:extremal} identify its limit with the polarized
orbifold at time $\tau_*$ of the same self-similar family. It is
polarized holomorphically isometric to the time-zero limit. The
configuration constructed at integer times therefore has the required
central polarized fibre for the original sequence as well.
\end{proof}

We use the convention that K-semistability requires $\DF\geq0$.
The following computation fixes the sign explicitly. For $A\coloneqq A_Z$, put
$\omega_A\coloneqq r\omega$, $S_A\coloneqq S/r$, and $dV_A\coloneqq r^2dV$.
Let $\eta$ be an algebraic $\C^*$-action on $Z$ linearized on $A$.
On $V_k\coloneqq H^0(Z,A^k)$, write the usual section action as
$\eta(e^s)=e^{sD_{\eta,k}}$, with integral weights
$\alpha_{k,1}(\eta),\ldots,\alpha_{k,d_k}(\eta)$. Define
\[
d_k\coloneqq\dim_\C V_k,\qquad
w_k(\eta)\coloneqq\tr D_{\eta,k}
=\sum_{j=1}^{d_k}\alpha_{k,j}(\eta).
\]
Thus $w_k(\eta)$ is the weight on $\det V_k$; unlike $d_k$, it
depends on the linearized action. For large $k$, write
\[
d_k=a_0k^2+a_1k+O(1),\qquad
w_k(\eta)=b_0(\eta)k^3+b_1(\eta)k^2+O(k),
\]
and define
\[
\DF_Z(\eta)\coloneqq\frac{a_1}{a_0}b_0(\eta)-b_1(\eta).
\]
This is the Donaldson--Futaki invariant of the product configuration
of $(Z,A)$ with action $a\cdot(q,z)=(\eta(a)q,az)$ and polarization
$\operatorname{pr}_Z^*A$. More generally, the invariant of a test
configuration equals $\DF_Z$ of its induced central-fibre action.
Dividing by $a_0$ is another convention and changes neither sign
nor vanishing.

Set $F_s\coloneqq\eta(e^{-s})$, using its given lift to $A$, and define
\[
\theta_A\coloneqq
-\left.\frac{d}{ds}\right|_{s=0}\log\frac{F_s^*h_A}{h_A}.
\]
Then
\[
\left.\frac{d}{ds}\right|_{s=0}F_s^*\omega_A=\ddc\theta_A,
\qquad
\left.\frac{d}{ds}\right|_{s=0}F_s^*dV_A
=(\Delta_A\theta_A)\,dV_A.
\]
Choose a complete orthonormal basis for the inner product defined by
$(h_A^k,dV_A)$, and write $G_k(s)$ for its Gram matrix with respect to
$F_s^*(h_A^k,dV_A)$. Pullback of sections by $F_s$ is the usual section
action of $\eta(e^s)$. The same change of variables as in
Lemma~\ref{tc:lem:gram} therefore gives
\[
\tr G_k'(0)
=-\tr D_{\eta,k}-\overline{\tr D_{\eta,k}}
=-2w_k(\eta).
\] Here the last equality uses the integral, hence real, weights of the
algebraic action; it does not require $D_{\eta,k}$ to be Hermitian.
Differentiating the defining integrals of the Gram matrix also gives
\[
2w_k(\eta)
=\int_Z(k\theta_A-\Delta_A\theta_A)\rho_k\,dV_A,
\]
where $\rho_k$ is the Bergman function of $(A^k,h_A^k,dV_A)$.

For every fixed smooth orbifold function $f$, the orbifold Bergman
expansion of \cite[Theorem 1.4 and (1.6), p.~4]{DLM} gives, in our
normalization,
\[
\int_Z f\rho_k\,dV_A
=(2\pi)^{-2}k^2\int_Z f\,dV_A
 +\frac12(2\pi)^{-2}k\int_Z fS_A\,dV_A+O_f(1).
\]
For comparison with that reference, rescale the K\"ahler form to
$\omega'=\omega_A/(2\pi)$. Its volume form is
$(2\pi)^{-2}dV_A$ and its Riemannian scalar curvature is
$4\pi S_A$, so the coefficient $r^X/(8\pi)$ there becomes $S_A/2$.
Indeed, the additional terms near the finitely many isolated quotient
points are bounded by
$Ck^2(1+\sqrt{k}\,d)^N e^{-c\sqrt{k}\,d}$, where $d$ is the distance
to the singular set. Their integrals in real dimension four are
$O(1)$. Thus they do not affect the two displayed coefficients.
Applying this formula to $f=\theta_A$ and to
$f=\Delta_A\theta_A$, and using orbifold integration by parts, yields
\[
\int_Z\Delta_A\theta_A\,dV_A=0,
\qquad
\int_Z(\Delta_A\theta_A)\rho_k\,dV_A=O(k).
\]
Comparison of coefficients in the exact trace identity, together with
the nonequivariant Riemann--Roch calculation above, proves
\begin{equation}\label{tc:eq:coefficients}
\begin{aligned}
a_0&=(2\pi)^{-2}\int_ZdV_A,&
a_1&=\tfrac12(2\pi)^{-2}\int_ZS_A\,dV_A,\\
b_0(\eta)&=\tfrac12(2\pi)^{-2}\int_Z\theta_A\,dV_A,&
b_1(\eta)&=\tfrac14(2\pi)^{-2}\int_Z\theta_A S_A\,dV_A.
\end{aligned}
\end{equation}
This argument applies even when the circle subgroup of $\eta$ does not
preserve the chosen metric.

With $\theta\coloneqq \theta_A/r$ this proves
\begin{equation}\label{tc:eq:futaki}
\DF_Z(\eta)=-\frac{r^2}{4(2\pi)^2}
\int_Z\theta(S-\bar S)\,dV.
\end{equation}
In particular $\DF_Z$ is linear on commuting generators and unchanged
by scalar linearization shifts. The zero-dimensional sector orders in
Lemma~\ref{tc:lem:spectrum} justify these formulas even for noncyclic local
quotient groups. The calculation uses the central orbifold and does not
require the total space to be $\Q$-Gorenstein.

The action $\lambda(e^{-s})=e^{bs\Lambda}$ on $Z$ is $F_{bs}$,
so increasing $s$ advances the limiting extremal flow.
Its lift potential therefore satisfies
\[
\ddc\theta_{\lambda,A}
=\left.\frac{d}{ds}\right|_0\lambda(e^{-s})^*\omega_A
=br\,\ddc S.
\]
Since pluriharmonic functions on the compact connected orbifold
are constant, $\theta_\lambda\coloneqq\theta_{\lambda,A}/r
=b(S-\bar S)+c$.
The constant reflects the bundle lift; trace-free matrix
normalization need not make the potential mean zero.

By Lemma~\ref{tc:lem:extremal}, the limiting average curvature is
$\bar S$. Put $\Ca(Z)\coloneqq\int_Z(S-\bar S)^2\,dV$.
The constant drops out of
$\int_Z\theta_\lambda(S-\bar S)\,dV=b\Ca(Z)$, so
\eqref{tc:eq:futaki} gives
\begin{equation}\label{tc:eq:negative}
\DF_Z(\lambda)=-c_b\Ca(Z),\qquad
c_b\coloneqq\frac{br^2}{4(2\pi)^2}>0.
\end{equation}
This is strictly negative for a non-cscK limit. The sign comes
from the forward-time choice in \eqref{tc:eq:lambda};
$b$ and $r$ are fixed before choosing $m$.

Both $\lambda$ and $\mu$ preserve the polarized central fibre
$(Z,A)$; fixing its Hilbert point does not mean fixing its points.
In particular, $\DF_Z(\lambda)$ concerns this action and does not
assert that $\lambda$ alone degenerates $M$ to $Z$.
Their commuting section generators give
\[
D_{\nu_m,k}=mD_{\lambda,k}+D_{\mu,k},\qquad
w_k(\nu_m)=m\,w_k(\lambda)+w_k(\mu).
\]
Comparing coefficients, with the same Hilbert coefficients
$a_0,a_1$ throughout, proves
\begin{equation}\label{tc:eq:dfsum}
\DF(\mathcal X_m,\mathcal L_m)
=\DF_Z(\nu_m)
=m\DF_Z(\lambda)+\DF_Z(\mu).
\end{equation}
If $Z$ is non-cscK, choose $m$ large enough both for the weight selection
in Lemma~\ref{tc:lem:combine} and for strict negativity in
\eqref{tc:eq:dfsum}. If $Z$ is cscK, every central-fibre Futaki character is
zero, so the direct configuration has $\DF=0$.

\begin{remark}[The two-stage description and the combined configuration]
\label{tc:rem:two-stage}
The construction first reaches $g[\overline M]$ by $\lambda$
and then $[Z]$ by $\mu$. If their commuting generators are
$A=-b\Lambda$ and $D$, the combined generator is $mA+D$,
or $A+D/m$ after positive rescaling. Thus the combined test configuration
has as its central fibre the polarized variety underlying the prescribed
geometric limit, but need not induce the same weighted filtration as the
Calabi flow. Although $D$ preserves the $A$-weight spaces, the weights
of $A+D/m$ may split them. We assert neither preservation of a normalized
invariant nor optimality. Neither $\mu$ nor $m$ is asserted to be canonical.

Chen--Sun--Wang likewise obtain a common generator and a second
degeneration in its centralizer for K\"ahler--Ricci soliton limits
\cite[Lemma 3.4 and Theorem 2.12]{CSW}. The first degeneration
records an intrinsic real filtration and its associated graded variety
\cite[Proposition 3.1 and the discussion after Proposition 3.6]{CSW}.
Soliton weights may be irrational; the rational approximation in
\cite[Section 3.1, p.~3159]{CSW} preserves the intermediate variety
but changes the grading. It can be chosen in the closure torus of
the soliton action and hence still commutes with the second generator.
Lemma~\ref{tc:lem:combine} then also gives an ordinary configuration
reaching the final variety, without asserting that it induces the
original weighted filtration or is optimal. In our polarized extremal setting,
rationality makes the first action algebraic without approximation;
compare \cite[Proposition 4.2, p.~3169]{CSW}.
\end{remark}

This proves part (2) of Theorem~\ref{app:thm:main}. We discuss the
consequence for stability after proving uniqueness of the flow limit.
The construction above depends on the complete-series identification in
Lemma~\ref{tc:lem:hilb-limit}; it does not require metric-limit uniqueness
or the global Bergman lower bound.

\subsection{Uniqueness of the flow limit}\label{app:sec:unique}\label{app:sec:fixed}

We now prove that all geometric limits of the flow are
polarized holomorphically isometric. We first recall the uniqueness
theorem for extremal metrics on a fixed orbifold, in a form that preserves
the common extremal field. Together with reductivity, this theorem shows
that the connected Hilbert limit set of
Subsection~\ref{app:sec:dynamics} consists of one compact-group orbit.
We then pass from integer-time limits to arbitrary real-time sequences
and deduce convergence in the geometric topology of the compactness
property.

We use Theorem~\ref{un:thm:fixed}, whose proof is given in
Subsection~\ref{app:sec:orbifold-uniqueness}. Its automorphism commutes
with the common extremal field and therefore preserves the matrix
centralizer used in Subsection~\ref{app:sec:dynamics}. The theorem
concerns a fixed orbifold; its proof requires no estimates uniform
over a degenerating family.

\begin{lemma}\label{un:lem:orbit}
Let $G$ be a complex reductive linear algebraic group, let $W$ be a
finite-dimensional complex representation of $G$, and let $K$ be a
maximal compact subgroup of $G$. Let $\mathscr C$ be a nonempty,
compact, connected, $K$-invariant subset of the closure of a single
$G$-orbit $O\subset\mathbb P(W)$. Suppose that, for every
$x\in\mathscr C$,
\begin{equation}\label{un:eq:orbit-hyp}
\Stab_G(x)\ \hbox{is reductive},\qquad
\mathscr C\cap Gx=Kx.
\end{equation}
Then $\mathscr C$ is one $K$-orbit.
\end{lemma}

\begin{proof}
We adapt the slice argument in
\cite[proof of Theorem 2.12, p.~3154]{CSW} to the intersection
hypothesis \eqref{un:eq:orbit-hyp}. We give the local closure argument
explicitly; the transverse construction below is the elementary one
in \cite[Section 5.2]{Donaldson} and does not require $Gx$ to be closed.

Choose $x\in\mathscr C$ with minimal stabilizer dimension, and put
$H\coloneqq\Stab_G(x)$ and $Y\coloneqq\overline O$, with its reduced
structure. Write $x=[v]$. Since $H$ is reductive, there is an
$H$-invariant splitting
\[
W=\C v\oplus W'.
\]
The corresponding $H$-invariant affine chart about $x$ is the
representation $\operatorname{Hom}(\C v,W')$, with $x$ corresponding
to the origin. Complete reducibility gives an $H$-invariant linear
complement to $T_x(Gx)$ in this chart. Let $S$ be the resulting
$H$-invariant affine linear subspace through $x$.

Consider
\[
F:G\times S\longrightarrow\mathbb P(W),\qquad F(g,s)=gs.
\]
Its differential is surjective at $(e,x)$. The rank condition and
$G$-equivariance therefore give an $H$-invariant Zariski-open
neighborhood $S^\circ\subset S$ of $x$ such that $F$ is smooth on
$G\times S^\circ$, of relative dimension $\dim H$. Set
\[
T\coloneqq Y\cap S^\circ,\qquad T_0\coloneqq O\cap S^\circ.
\]
Since $O$ is open and dense in $Y$ and smooth morphisms are open,
\[
F^{-1}(O)=G\times T_0
\quad\hbox{is dense in}\quad
F^{-1}(Y)=G\times T.
\]
Thus $T_0$ is dense in $T$. Smooth base change over the smooth orbit
$O$ also shows that $T_0$ is smooth, of dimension
\[
\dim T_0=\dim O-\dim G+\dim H.
\]

For $s\in T_0$, the kernel of $dF_{(e,s)}$ has dimension $\dim H$
and contains the subspace
\[
\{(\eta,-\eta\cdot s):\eta\in\Lie H\}.
\]
The two subspaces are therefore equal. Applying this to
$(\xi,0)$ for $\xi\in\Lie\Stab_G(s)$ gives
\[
\Lie\Stab_G(s)=\Lie\Stab_H(s).
\]
Consequently
\[
\dim(Hs)=\dim H-\dim\Stab_G(s)
       =\dim O-\dim G+\dim H
       =\dim T_0.
\]
Every $H^0$-orbit in $T_0$ is thus open. Since the finite-type smooth
variety $T_0$ has finitely many irreducible components, and an
irreducible component cannot contain two disjoint nonempty open
orbits, $T_0$ is a finite union of $H$-orbits.

By density, $T$ is the union of the closures in $S^\circ$ of these
finitely many orbits. Shrink about $x$ to exclude the closures that
do not contain $x$. For every remaining slice point $s$, all
$H$-invariant polynomials on the affine space $S$ take their value
at $x$: the point $s$ lies in the closure of an $H$-orbit whose
closure contains $x$. Therefore $s$ lies in the fibre of $x$ for the
affine quotient $S\to S/\!/H$. Each fibre of an affine reductive
quotient contains a unique closed orbit, which is contained in the
closure of every orbit in that fibre. Here $x$ is $H$-fixed, so that
closed orbit is $\{x\}$. Hence $x\in\overline{Hs}$ for every such $s$.

Finally, the submersion $F$ at $(e,x)$ has a local holomorphic
section. It follows that every $y\in Y$ sufficiently near $x$ can be
written as $y=gs$, with $g$ near $e$ and $s\in T$ in the neighborhood
just chosen. The $G$-invariance of $Y$ ensures that this slice point
belongs to $T$. We have proved
\begin{equation}\label{un:eq:closure-order}
x\in\overline{Gy}
\qquad\hbox{for every $y\in Y$ sufficiently near $x$}.
\end{equation}

If $y\in\mathscr C$ is sufficiently near $x$ and $Gy\ne Gx$, the
boundary-orbit dimension inequality gives
\[
\dim Gy>\dim Gx,
\qquad
\dim\Stab_G(y)<\dim\Stab_G(x),
\]
contrary to the choice of $x$. Thus nearby points of $\mathscr C$
lie in $Gx$, and \eqref{un:eq:orbit-hyp} places them in $Kx$.
By $K$-invariance, translates of this neighborhood make $Kx$ open in
$\mathscr C$. Compactness of $K$ makes $Kx$ closed. Connectedness
therefore gives $\mathscr C=Kx$.
\end{proof}

\begin{remark}
The necessary condition is the intersection formula
$\mathscr C\cap Gx=Kx$. It does not assert $Gx=Kx$ in the entire
ambient projective space. In the slice argument the relevant closure
is a complex-group orbit closure, not a compact-group orbit closure.
These distinctions also match the geometric verification in
\cite[Section 3.2, p.~3163; Section 4]{CSW}. Lemma~\ref{un:lem:orbit}
provides the precise version used here, independently of abbreviated
orbit notation in that source.
\end{remark}

We now verify each hypothesis of Lemma~\ref{un:lem:orbit} for the
limits of the fixed Calabi flow.

\textit{Step 1: the groups and the limit set.}
For the eigenspace decomposition $\C^N=\bigoplus_\alpha U_\alpha$
of $\Lambda$, the groups
\[
G_\Lambda=\prod_\alpha GL(U_\alpha),\qquad
K_\Lambda=\prod_\alpha U(U_\alpha)
\]
are reductive and maximal compact, respectively. They act on the
fixed projective representation of the Hilbert scheme.
Lemma~\ref{app:lem:limitset} gives a nonempty compact connected
$K_\Lambda$-invariant set
$\mathscr C\subset\overline{G_\Lambda[\overline M]}$.
Each point represents an extremal orbifold embedded by a complete
orthonormal basis, with $\Lambda$ inducing its extremal flow.
In particular, connectedness comes from the compatible dynamics,
before any use of uniqueness.

\textit{Step 2: reductivity of the stabilizers.}
For $[Z]\in\mathscr C$, write
$H_Z\coloneqq\Stab_{G_\Lambda}([Z])$ and let $Y_Z$ be the
holomorphic extremal field. The complete embedding identifies the
projective stabilizer with $\Aut(Z,A_Z)$. Here the orbifold structure
is determined by the coarse space: a local isomorphism lifts on the
punctured quotient charts to their simply connected finite covers,
and Hartogs extends the lift and its inverse across the origin.
The local groups act freely away from the origin, so these lifts
identify the orbifold charts.

An element of $H_Z$ induces a polarized automorphism commuting
with $Y_Z$. Conversely, such an automorphism has a linear
representative $B$ in the complete embedding, and commutation with
the extremal flow gives
$B\Lambda B^{-1}=\Lambda+cI$. Indeed, a linear endomorphism whose
induced projective vector field vanishes on the connected,
nondegenerate embedded $Z$ must be scalar.
Taking traces gives $c=0$.
Thus there is an exact sequence
\[
1\longrightarrow\C^*\longrightarrow H_Z
\longrightarrow\Aut(Z,A_Z,Y_Z)\longrightarrow1,
\]
where the kernel consists of scalar matrices.
Lemma~\ref{tc:lem:red} makes the quotient reductive. The unipotent
radical of $H_Z$ therefore lies in this scalar torus and is trivial.
Hence every stabilizer $H_Z$ is reductive, as required.

\textit{Step 3: the orbit-intersection condition.}
Fix $[Z_1]\in\mathscr C$. The inclusion
$K_\Lambda[Z_1]\subset\mathscr C\cap G_\Lambda[Z_1]$
follows from $K_\Lambda$-invariance. For the reverse inclusion,
suppose $[Z_2]=g[Z_1]$ with $g\in G_\Lambda$ and $[Z_2]\in\mathscr C$.
Choose the associated extremal metrics $\omega_i$ and Hermitian
metrics $h_i$ on $A_i$, with $\sqrt{-1}F_{h_i}=r\omega_i$.
The projective transformation restricts to a polarized orbifold
isomorphism $g:Z_1\to Z_2$, by the lifting argument in Step~2.
It need not be an isometry.

Since $g\Lambda=\Lambda g$, it intertwines the two extremal fields.
Consequently $g^*\omega_2$ and $\omega_1$ have the same holomorphic
extremal field on $Z_1$. Moreover,
\[
[g^*\omega_2]=[\omega_1]=\frac{2\pi}{r}c_1(A_1).
\]
The remaining numerical hypothesis $q(Z_1)=0$ of
Theorem~\ref{un:thm:fixed} follows from
Lemma~\ref{tc:lem:hilbert}. That theorem supplies
$a\in\Aut^0(Z_1,A_1,Y_1)$ such that
\[
a^*g^*\omega_2=\omega_1.
\]
Thus $F\coloneqq g\circ a$ is a polarized holomorphic isometry
which still intertwines the extremal fields.

Under a bundle isomorphism $F^*A_2\simeq A_1$, the metrics
$F^*h_2$ and $h_1$ have the same curvature. Their logarithmic ratio
is a constant, so rescaling the lift makes it an isometry.
Since $F$ also preserves volume, its pullback on complete section
spaces is unitary. The dual map therefore has a unitary matrix $k$
in the chosen orthonormal coordinates, with $k[Z_1]=[Z_2]$.

Finally, $F$ intertwines the projective extremal flows. By the
complete embedding, their matrices satisfy
$k\Lambda k^{-1}=\Lambda+cI$; taking traces again gives $c=0$.
Thus $k\in U(N)\cap G_\Lambda=K_\Lambda$. This proves
\[
\mathscr C\cap G_\Lambda[Z_1]=K_\Lambda[Z_1].
\]
Notice that the unitary representative is obtained from the
corrected map $F=g\circ a$; the original map $g$ need not be unitary.

\textit{Step 4: the geometric conclusion.}
Lemma~\ref{un:lem:orbit} now gives
$\mathscr C=K_\Lambda[Z_\infty]$. The correction in Step~3
identifies the metrics, including two metrics represented by the
same Hilbert point (take $g=I$). Thus all integer-time limits are
polarized holomorphically isometric.

For arbitrary $t_j\to\infty$, write $t_j=n_j+\tau_j$ and extract
$\tau_j\to\tau\in[0,1]$. Compactness at $n_j$ on common time
windows identifies the limit with $F_\tau^*\omega_\infty$, in the
same polarized isometry class. Relative sequential compactness and
uniqueness then imply convergence in the quotient geometric topology:
otherwise a sequence outside a fixed neighborhood would have a
subsequence converging to that class. An isometry of extremal
metrics identifies their scalar curvatures and conjugates their
flows $F_\tau$, so it identifies the entire limiting families.
This proves part (3) of Theorem~\ref{app:thm:main}, including its
translated-window formulation; it does not select a fixed global
coordinate gauge.

\subsection{Stability and existence}\label{app:sec:stability}

We conclude by relating the geometric limit of the Calabi flow to
the stability of the original polarized surface. The test configuration
constructed in Subsection~\ref{app:sec:tc} shows that K-semistability forces
the limit to have constant scalar curvature. K-polystability further
identifies the limit with the original surface and removes its orbifold
singularities.

\begin{corollary}[Stability and the smooth limit]\label{app:cor:stable}
Under the assumptions of Theorem~\ref{app:thm:main}, if $(M,L)$ is
K-semistable, its unique limiting orbifold metric satisfies
$S_{\omega_\infty}\equiv\bar S$.
If $(M,L)$ is K-polystable, there is a polarized isomorphism
$F:(M,L^r)\to(Z_\infty,A_\infty)$, and $F^*\omega_\infty$ is a smooth
cscK metric on the original complex surface in $2\pi c_1(L)$.
The Calabi flow converges to its polarized holomorphic-isometry class
in the geometric sense of part~(3) of the theorem.
\end{corollary}
\begin{proof}
A non-cscK limit gives a test configuration with negative DF by
Subsection~\ref{app:sec:tc}, contradicting K-semistability. For a cscK limit,
the configuration has DF zero. K-polystability forces a product, so
$(Z,A_Z)\simeq(M,L^r)$. In particular, the coarse complex space $Z$
is smooth.

We show directly that every local orbifold group is trivial.
Let $q\in Z$ have a uniformizing chart $B/\Gamma$, where
$B\subset\C^2$ is a $\Gamma$-invariant ball centered at $0$ and
$\Gamma\subset U(2)$ acts freely on $B\setminus\{0\}$, as required
by the compactness property. Shrink $B$ so that $B/\Gamma$ lies
in a holomorphic coordinate chart of the smooth surface $Z$.
Composing the quotient map with these coordinates gives a
holomorphic map $\Psi:B\to\C^2$ satisfying
$\Psi\circ\gamma=\Psi$ for every $\gamma\in\Gamma$.
Freeness away from $0$ makes $\Psi$ locally biholomorphic there,
so its Jacobian determinant
$J_\Psi\coloneqq\det D\Psi$ is nonzero on $B\setminus\{0\}$.
By Hartogs' theorem, $1/J_\Psi$ extends holomorphically to a
function $f$ on $B$. The identity $J_\Psi f=1$ then holds on
all of $B$, so $D\Psi(0)$ is invertible.
Differentiating $\Psi\circ\gamma=\Psi$ at $0$ gives
\[
D\Psi(0)\,\gamma=D\Psi(0),
\]
and hence $\gamma=I$ for every $\gamma\in\Gamma$.
Thus all local groups are trivial. The limiting metric, already
smooth in orbifold charts, is therefore a smooth metric on $Z$.
Its pullback is a smooth cscK metric in the original class. The uniqueness
and convergence conclusion is part (3) already proved.
\end{proof}
This proves part (4) and completes Theorem~\ref{app:thm:main}.

\section{Finite-group equivariant applications}\label{eqv:sec:main}

We prove finite-group equivariant versions of the Bergman estimate,
test-configuration construction and limit uniqueness. They yield a
cscK existence criterion testing only equivariant configurations
with integral central fibre.

\begin{lemma}[Linearization and rescaling]\label{eqv:lem:linearization}
Let a finite group $\Gamma$ act holomorphically on the connected
polarized surface $(M,L)$, with $\gamma^*L\simeq L$ for every
$\gamma\in\Gamma$. Suppose that $\omega(t)\in2\pi c_1(L)$ is a
global smooth Calabi flow with $\Gamma$-invariant initial metric.
Put $d\coloneqq|\Gamma|$, or $d=1$ if $L$ is already linearized.
Then $L^{\otimes d}$ admits a genuine $\Gamma$-linearization: a
holomorphic fibrewise linear action covering the given action on $M$
and satisfying the group law. The rescaled data
\begin{equation}\label{eqv:eq:rescaling}
\widetilde L\coloneqq L^{\otimes d},\qquad
\widetilde\omega(\tau)\coloneqq d\,\omega(\tau/d^2),\qquad
\widetilde h(\tau)\coloneqq h(\tau/d^2)^{\otimes d}
\end{equation}
form a $\Gamma$-invariant polarized Calabi flow, where
$\sqrt{-1}F_{h(t)}=\omega(t)$.
The Sobolev bound \eqref{tc:eq:sob} and numerical condition
\eqref{tc:eq:numerical} are preserved. A cscK metric for
$\widetilde L$, divided by $d$, lies in the original class;
a test configuration of exponent $r$ for $\widetilde L$ has
exponent $dr$ for $L$.
\end{lemma}
\begin{proof}
For $d=|\Gamma|$, reindexing the factors of
$N\coloneqq\bigotimes_{\gamma\in\Gamma}\gamma^*L\simeq L^{\otimes d}$
satisfies the group law and gives a canonical $\Gamma$-linearization.
Transfer it through a fixed isomorphism to $L^{\otimes d}$;
for $d=1$, use the given linearization. Uniqueness of smooth Calabi flow preserves
$\Gamma$-invariance of $\omega(t)$. For each $\gamma$, the ratio
$\gamma^*\widetilde h/\widetilde h$ has pluriharmonic logarithm,
hence is a positive constant; finite order of $\gamma$ makes it one.
Since $S_{d\omega}=d^{-1}S_\omega$,
\[
\partial_\tau\widetilde\omega
=d^{-1}\partial_t\omega
=d^{-1}\ddc S_\omega
=\ddc S_{\widetilde\omega}.
\]
Tensor powers give the curvature normalization. The scale-invariant
Sobolev constant and $\chi(\mathcal O_M)$ are unchanged, while
$c_1(M)\cdot c_1(L)$ is multiplied by $d$. The last assertions
follow directly from the rescaling and the definition of exponent.
\end{proof}

By Lemma~\ref{eqv:lem:linearization}, we may assume that $L$ has a
fixed $\Gamma$-linearization preserving $h(t)$, and retain the notation
$L,\omega(t),h(t)$ for the rescaled data.
A $\Gamma$-equivariant test configuration means that the action lifts
to its polarization, acts trivially on the base, commutes with the
$\C^*$ action, and restricts over $\C^*$ to the prescribed action
on the polarized general fibre.

\begin{theorem}[Equivariant applications]\label{eqv:thm:main}
Let $(M,L)$ be a smooth connected polarized complex surface with a
fixed linearized action of a finite group $\Gamma$. Let
$\omega(t)\in2\pi c_1(L)$, $t\geq0$, be a global smooth Calabi flow
satisfying \eqref{tc:eq:sob} and \eqref{tc:eq:numerical}, with
$\Gamma$-invariant initial metric. Then the entire flow is
$\Gamma$-invariant. There is a common integer $r>0$, independent of
the subsequence defining a geometric limit, with the following properties.
\begin{enumerate}[font=\normalfont\bfseries]
\item \textbf{\upshape Bergman estimate and equivariant embeddings.}
There are $0<b\leq B<\infty$ such that
$b\leq\rho_{r,t}(x)\leq B$ for all $t\geq0$ and $x\in M$.
For every equivariant geometric limit $(Z,\omega_Z)$, the bundle
$L_Z^r$ descends to a very ample Cartier line bundle $A_Z$, with
$[\omega_Z]=(2\pi/r)c_1(A_Z)$ and an induced $\Gamma$-linearization.
Along any geometrically convergent subsequence, after further
extraction, complete orthonormal bases converge smoothly on regular
compact subsets to a complete limit basis. The complete-section
identifications intertwine the $\Gamma$ representations.
The resulting embeddings of $M$ and $Z$
into $\mathbb P(H^0(M,L^r)^*)$ are equivariant for one fixed
representation. The estimate uses the full section space.
If the action on $M$ is effective, so is the limiting action on $Z$.
\item \textbf{\upshape Equivariant test configurations.}
Every $(Z,A_Z)$ in part~(1), with its limiting group action, is the
central polarized fibre of a $\Gamma$-equivariant normal, relatively
ample test configuration for $(M,L)$ of exponent $r$.
The entire central fibre is a normal, reduced and irreducible scheme.
The configuration can be chosen with strictly negative
Donaldson--Futaki invariant if $\omega_Z$ is not cscK, and with zero
invariant if $\omega_Z$ is cscK.
\item \textbf{\upshape Equivariant uniqueness and convergence.}
For any two equivariant geometric limits $(Z_i,A_i,\omega_i)$,
with actions denoted by $\gamma_i$, there is an orbifold
biholomorphism $F:Z_1\to Z_2$ such that
\[
F\circ\gamma_1=\gamma_2\circ F\quad(\gamma\in\Gamma),\qquad
F^*A_2\simeq A_1,\qquad F^*\omega_2=\omega_1.
\]
This fixed flow converges to the resulting equivariant polarized
holomorphic-isometry class, in particular in the equivariant
Gromov--Hausdorff sense. Every sequence $t_j\to\infty$ has a
subsequence with the coherent regular smooth convergence on bounded
translated time intervals of Theorem~\ref{app:thm:compactness} to
the self-similar flow through the unique limit. The GH correspondences
are equivariant, and the local group-action maps converge smoothly
as in Lemma~\ref{eqv:lem:compactness}.
\end{enumerate}
\end{theorem}
We shall apply Lemma~\ref{eqv:lem:compactness} to the invariant flow and its
time translates. The limiting action preserves scalar curvature and
its gradient, hence commutes with the extremal transport $F_\tau$.
We perform the polarized constructions in the representation centralizer. 

The invariant Hermitian metrics from Lemma~\ref{eqv:lem:linearization}
give invariant $L^2$ products $H_t$ on $H^0(M,L^r)$.
In the fixed dual coordinates determined by an $H_0$-orthonormal
ordered basis, write
\[
\mathbb V\coloneqq H^0(M,L^r)^*
 =\bigoplus_\rho W_\rho\otimes\C^{m_\rho}
\]
for the fixed evaluation-coordinate representation of $\Gamma$.
Here $W_\rho$ are its distinct irreducible representations.
In the induced $H_0$-dual unitary coordinates put
\begin{equation}\label{eqv:eq:centralizer}
\mathcal R\coloneqq \operatorname{Cent}_{\operatorname{GL}(\mathbb V)}(\Gamma)
 =\prod_\rho\operatorname{GL}(m_\rho),\qquad
K_{\mathcal R}\coloneqq \mathcal R\cap U(\mathbb V)=\prod_\rho U(m_\rho).
\end{equation}
All Gram matrices, their positive square roots and the successive
orthonormalizations in \eqref{tc:eq:recursion} lie in $\mathcal R$.
Indeed, if $U_\gamma$ is the unitary evaluation-coordinate action,
invariance gives $U_\gamma^{-1}G_t(U_\gamma^{-1})^*=G_t$, so
$G_t$ commutes with $U_\gamma$.
The global determinant normalization is scalar and preserves this fact.

\begin{corollary}[Equivariant polarization and complete sections]
\label{eqv:cor:polarization}
At the common power $r$ of Lemma~\ref{tc:lem:hilbert}, the limiting
polarization $A_Z\coloneqq L_Z^r$ has a genuine $\Gamma$-linearization covering
the action in Lemma~\ref{eqv:lem:compactness} and preserving its
Hermitian metric. The complete section identifications intertwine the
fixed $\Gamma$ representation, and the limiting projective embeddings
are equivariant. If the action on $M$ is effective, so is its action
on $Z$.
\end{corollary}
\begin{proof}
Lemma~\ref{tc:lem:sections} applies to every ordered orthonormal basis,
so we may use the bases obtained by the centralizer Gram changes
above. The source complete embeddings
$\iota_j:M\to\mathbb P(\mathbb V)$ are equivariant for one fixed
projective representation. On regular compact subsets, complete-section
convergence and Lemma~\ref{eqv:lem:compactness} pass their identities to
\[
\iota_\infty\circ\gamma_\infty
 =\gamma\circ\iota_\infty.
\]
Projective evaluation removes the local scalar ambiguity of the bundle
gauges. Density of the regular locus extends the identity over $Z$.
The canonical linearization of $\mathcal O(1)$, restricted to this
embedding, gives the required genuine action on
$A_Z=\iota_\infty^*\mathcal O(1)$, and hence the stated intertwining
of complete section spaces.

Its curvature is $r\omega_Z$, which the action preserves. The ratio
of the pulled-back Hermitian metric to the original one therefore has
pluriharmonic logarithm and is a positive constant on the compact
connected orbifold. These constants form a character
$\Gamma\to\R_{>0}$, necessarily trivial. Thus the lifted action is
unitary. Finally, if $\gamma_\infty$ is the identity, the fixed
projective transformation fixes $\iota_\infty(Z)$ pointwise. This
irreducible nondegenerate image must lie in one projectivized
eigenspace of the finite-order transformation; nondegeneracy makes
that eigenspace all of $\mathbb V$. Its linear representative is
scalar, so its projective action also fixes the source embedding pointwise.
Effectiveness on $M$ then gives $\gamma=1$.
\end{proof}
The divisor argument in Subsection~\ref{app:sec:partial} still uses the
full section space $\mathbb V^*$; invariant sections alone need not
define an embedding.

\begin{lemma}[Common spectra in the multiplicity spaces]\label{eqv:lem:spectra}
There is a globally trace-free Hermitian
$\Lambda\in\operatorname{Lie}(\mathcal R)$ such that every limiting
Gram generator is $K_{\mathcal R}$-conjugate to $\Lambda$ after removing
its global trace. Some positive multiple of $\Lambda$ has integral
eigenvalues, and \eqref{tc:eq:star} holds with conjugating unitaries in
$K_{\mathcal R}$.
\end{lemma}
\begin{proof}
Lemma~\ref{tc:lem:gram} gives the bounded-window limit $e^{-2\tau C_Z}$,
and invariance implies
$C_Z=\bigoplus_\rho I_{W_\rho}\otimes C_{Z,\rho}$.
The rationality argument of Lemma~\ref{tc:lem:spectrum} applies to its
globally trace-free eigenvalues with the same fixed normalization.

For real $t$, choose the canonical positive orthonormalization $s_t$ of
$H_t$ and set
\[
D(t)\coloneqq -\tfrac12\log\operatorname{Gram}_{H_{t+1}}(s_t),\qquad
D^0(t)\coloneqq D(t)-\frac{\tr_{\mathbb V}D(t)}{\dim\mathbb V}I.
\]
Record the ordered spectrum of $D^0(t)$ on each multiplicity space.
The resulting vector is continuous and has relatively compact tails
by Lemma~\ref{tc:lem:gram}. Its cluster set is connected. Every cluster
value is the collection of block spectra of a globally trace-free
limiting generator, so Lemma~\ref{tc:lem:spectrum} places it in
$\Q^{\sum_\rho m_\rho}$ with the normalization already fixed above.
The cluster set therefore consists of a single point.
This determines the common block spectra and the required
$K_{\mathcal R}$-conjugacy. Clearing the finitely many rational
denominators gives a positive integer $b$ such that $b\Lambda$ has
integral eigenvalues. The Gram limits at times $1,2$, followed by
reorthonormalization at time $1$, give the two-step condition exactly
as in Subsection~\ref{app:sec:dynamics}.
\end{proof}
It is essential to remove only the global trace
$\sum_\rho(\dim W_\rho)\tr C_{Z,\rho}$. Removing a separate trace on
each block would change the relative projective weights. Conversely,
knowing only the total spectrum would not prevent eigenvalues from
moving between inequivalent representation blocks.

Set
\[
\mathcal R_\Lambda\coloneqq \operatorname{Cent}_{\mathcal R}(\Lambda),\qquad
K_{\mathcal R,\Lambda}\coloneqq K_{\mathcal R}\cap\mathcal R_\Lambda.
\]
The latter is a connected product of unitary groups on the joint
multiplicity and eigenvalue spaces.

\begin{lemma}\label{eqv:lem:orbit}
The construction of Lemma~\ref{app:lem:limitset} can be carried out
inside $\mathcal R$. In particular, one can choose $C_0\in\mathcal R$
so that the compatible Hilbert limit set $\mathscr C_\Gamma$ is compact,
connected, $K_{\mathcal R,\Lambda}$-invariant, and satisfies
\[
[\bar M]=\lim_{t\to\infty}e^{t\Lambda}C_0[M],\qquad
\mathscr C_\Gamma\subset
\overline{\mathcal R_\Lambda[\bar M]}.
\]
\end{lemma}
\begin{proof}
We verify that the coordinate choices in \cite[Section 2]{CSW} preserve
the representation blocks. The gauge construction uses compactness,
multiplication, and paths in the unitary centralizer; these remain in
$K_{\mathcal R}$ and $K_{\mathcal R,\Lambda}$. The growth filtrations
are defined using the $\widehat A_i$, those gauges and sums of $\Lambda$
eigenspaces, so they and their Grassmannian limits are
$\Gamma$-invariant. Their complementary splitting spaces $W_s$ are
invariant as well.

Write $U_s$ for the corresponding $\Lambda$ eigenspaces. For the
chosen gauges $g_i\in K_{\mathcal R}$, the construction gives
$g_i^{-1}\widehat A_iW_s\to U_s$. For one sufficiently large $i$, the maps
\[
P_{U_s}g_i^{-1}\widehat A_i|_{W_s}:W_s\longrightarrow U_s
\]
are invertible intertwiners: orthogonal projection between close
invariant subspaces is an invertible intertwiner. Their direct sum gives
$C_0\in\mathcal R$. The near-identity corrections obtained by projection
on the converging gauged splittings also lie in $\mathcal R$.
The fixed $W_s$ themselves need not be close to $U_s$.
After correction the matrices commute with $\Lambda$ as well.
The exterior-power growth and weight projection argument of
\cite[Proposition 2.10]{CSW} therefore gives the orbit inclusion in
$\mathcal R_\Lambda$. The connectedness argument of its Lemma 2.6
uses the connected group $K_{\mathcal R,\Lambda}$ and gives the
asserted connected limit set. Thus we have checked the construction for
the specific product of general linear groups in
\eqref{eqv:eq:centralizer}.
\end{proof}

Choose $b>0$ with $b\Lambda$ integral and put
$\lambda(z)\coloneqq \exp(-b\Lambda\log z)\in\mathcal R$.
The stabilizer of $[Z]$ in $\mathcal R_\Lambda$ is reductive.
Indeed, the extremal centralizer is the complexification of its compact
Hamiltonian isometry group, by Lemma~\ref{tc:lem:red}. Conjugation by
$\Gamma$ preserves that compact real form, since $\Gamma$ acts by
isometries. Taking its fixed Lie algebra preserves the complexification
property. The linear stabilizer adds a central scalar group; possible
differences between linear and projective commuting are finite
characters of $\Gamma$ and do not affect reductivity.

Reductive-stabilizer accessibility inside $\mathcal R_\Lambda$ now
gives $g\in\mathcal R_\Lambda$ and a one-parameter subgroup
$\mu:\C^*\to\mathcal R_\Lambda$ with
$\lim_{z\to0}\mu(z)g[\bar M]=[Z]$. Lemma~\ref{tc:lem:combine}
combines the commuting weights:
\[
\nu_m\coloneqq \lambda^m\mu\in\mathcal R,\qquad
\lim_{z\to0}\nu_m(z)gC_0[M]=[Z]
\]
for sufficiently large integers $m$. If $\Lambda=0$, use accessibility
directly on the closure of $\mathcal R[M]$.
The universal Hilbert family is $\Gamma$-equivariant because both its
initial coordinate change and its one-parameter subgroup commute with
$\Gamma$. The complete-series identification and normality proof in
Subsection~\ref{app:sec:tc} give the whole normal integral central scheme
$Z$ and a normal relatively ample total family.
The same Futaki calculation gives
\[
\DF_Z(\lambda)=-c_b\Ca(Z),\qquad
\DF(\mathcal X_m,\mathcal L_m)
=m\DF_Z(\lambda)+\DF_Z(\mu).
\]
This proves the asserted sign. Comparison of real times with their
integer parts uses the limiting extremal motion, which commutes with
$\Gamma$, and therefore preserves equivariance.

We now apply the finite-group uniqueness lemma from
Subsection~\ref{app:sec:orbifold-uniqueness}.
If $[Z_2]=g[Z_1]$ with $g\in\mathcal R_\Lambda$, the metrics
$\omega_1,g^*\omega_2$ have the same extremal field and are
$\Gamma$-invariant. Lemma~\ref{eqv:lem:fixed} gives the correcting
automorphism $a$ centralizing both $\Gamma$ and the field. Lift its path
from the identity to the complete section space. Each commutator with
the genuine $\Gamma$ representation is a scalar character of the finite
group, hence is identically one along this path. The polarized isometry
$g\circ a$ can be lifted unitarily, as in Subsection~\ref{app:sec:unique}.
The common global trace normalization makes its matrix commute with
$\Lambda$. We have proved
\[
\mathscr C_\Gamma\cap\mathcal R_\Lambda[Z]
=K_{\mathcal R,\Lambda}[Z].
\]
Lemma~\ref{un:lem:orbit} now gives one compact-group orbit. The same
integer-to-real-time comparison proves equivariant uniqueness on all
translated windows and full convergence in the stated quotient topology.
This completes the proof of Theorem~\ref{eqv:thm:main}.

\begin{corollary}\label{eqv:cor:restricted}
Under Theorem~\ref{eqv:thm:main}, suppose DF is nonnegative on every
$\Gamma$-equivariant normal ample test configuration whose central
fibre is integral. Then the unique limit is cscK. If in addition every
zero-DF configuration in this restricted class is a product, the
original polarized surface admits a smooth cscK metric in the given
class. If these products are equivariantly products, the identification
with the original surface can be chosen $\Gamma$-equivariantly. This
last assertion holds in particular when $\Aut^0(M)=\{1\}$.
\end{corollary}
\begin{proof}
Theorem~\ref{eqv:thm:main} produces precisely a configuration of the
type being tested. A negative DF excludes a non-cscK limit. For a cscK
limit its DF is zero, so a product identifies the polarized coarse
space with $M$. The Hartogs argument in
Corollary~\ref{app:cor:stable} removes the orbifold singularities and
gives the smooth cscK metric. An equivariant product retains the group
action. If $\Aut^0(M)=\{1\}$, any product $\C^*$ action on $M$ is
trivial; in a product trivialization the commuting fibrewise
$\Gamma$ action is consequently independent of the base coordinate.
Thus the identification is equivariant in this case as well.
\end{proof}
The criterion requires no uniform stability estimate and does not
assert full K-polystability.

\section{Examples of smooth cscK and extremal limits}
\label{ex:sec:examples}

We give two nonstationary examples satisfying the hypotheses of
Theorem~\ref{app:thm:main}. Their limiting metrics are already known.
The first limit is cscK but not K\"ahler--Einstein; the second has
nonconstant scalar curvature. All energy comparisons use the complex
scalar curvature and the volume form in \eqref{app:eq:conventions}.

For an ample class $A=c_1(L)$, put
\begin{equation}\label{ex:eq:general-threshold}
 c\coloneqq c_1(X),\qquad
 \delta(A)\coloneqq c^2-\frac{2(c\cdot A)^2}{3A^2}.
\end{equation}
When $c\cdot A>0$ and $\delta(A)>0$, the sufficient initial-energy
condition in \cite[Corollary 1.3]{LongTime} is
\begin{equation}\label{ex:eq:general-startup}
 \Ca(\omega_0)<8\pi^2\delta(A),\qquad [\omega_0]=2\pi A.
\end{equation}
It supplies global existence and a uniform Sobolev bound. Both surfaces
below are rational, so $\chi(\mathcal O_X)=1$.
For comparison with that reference, its volume is $\int_X\omega^2=2V$
and its normalized Calabi energy is $\Ca/V$, where
$V=\int_XdV$. Its bound $\Ca/V<16\pi^2\delta(A)/(2V)$
is exactly \eqref{ex:eq:general-startup}.

\subsection{A product cscK limit}

Let $X=\mathbb P^1\times\mathbb P^1$ and
$L=\mathcal O(1,2)$. Write $A_1,A_2$ for the pullbacks of the
hyperplane class on the two factors. Then
\[
 A_1^2=A_2^2=0,\quad A_1\cdot A_2=1,\quad
 c=2A_1+2A_2,\quad A=A_1+2A_2.
\]
In particular,
\begin{equation}\label{ex:eq:product-threshold}
 c^2=8,\qquad A^2=4,\qquad c\cdot A=6,\qquad
 \delta(A)=2.
\end{equation}
Normalize the Fubini--Study form $\omega_{\rm FS}$ on $\mathbb P^1$
by $[\omega_{\rm FS}]=2\pi c_1(\mathcal O(1))$. Its complex scalar
curvature is $2$. With pullbacks understood, the product metric
\[
 \omega_*=\omega_{{\rm FS},1}+2\omega_{{\rm FS},2}\in2\pi A
\]
has scalar curvature $2+1=3$. It is cscK but not K\"ahler--Einstein,
since its two factors have different Einstein constants.

Choose a smooth real function $f$ outside the kernel of the
Lichnerowicz operator $\mathcal L_*$ at $\omega_*$. For example,
one may take the product of the zero-mean height Hamiltonians on
the two spheres. For sufficiently small $\varepsilon$, set
\[
 \omega_\varepsilon=\omega_*+\varepsilon\ddc f.
\]
This form is positive and belongs to the same integral class.
The scalar-curvature linearization and the fixed average give
\begin{equation}\label{ex:eq:product-energy}
 \Ca(\omega_\varepsilon)
 =\varepsilon^2\int_X(\mathcal L_*f)^2\,dV_*
   +O(\varepsilon^3).
\end{equation}
The leading coefficient is positive. Consequently every sufficiently
small nonzero $\varepsilon$ gives
\[
 0<\Ca(\omega_\varepsilon)<16\pi^2
       =8\pi^2\delta(A).
\]
The flow is nonstationary and satisfies the global existence and
Sobolev hypotheses of Theorem~\ref{app:thm:main}.
Chen--He stability \cite[Theorem 1.2 and Remark 1.3]{CHFlow}
also gives smooth
convergence to a nearby cscK metric. Fixed-class uniqueness identifies
its polarized holomorphic-isometry class with that of $\omega_*$.
Indeed, an automorphism preserving $A_1+2A_2$ cannot interchange
the factors and preserves $\mathcal O(1,2)$.
Thus the unique geometric limit in Theorem~\ref{app:thm:main} is the
smooth product cscK metric, up to this equivalence.

\subsection{A non-cscK extremal limit on the first Hirzebruch surface}

Let $X=\mathbb F_1=\operatorname{Bl}_p\mathbb P^2$ and
$L=-K_X$. If $H$ is the pullback of the line class and $E$ is the
exceptional class, then
\begin{equation}\label{ex:eq:f1-threshold}
 A=c=3H-E,\qquad A^2=c\cdot A=c^2=8,\qquad
 \delta(A)=\frac83.
\end{equation}
Calabi's construction gives a $U(2)$-invariant extremal metric
$\omega_{\rm ext}\in2\pi A$; see \cite[Section 7]{LeBrunToric}.
We compute its energy to check \eqref{ex:eq:general-startup} directly.

Choose moment coordinates with angular periods $2\pi$ and polygon
\[
 P=\{(x,y):x\geq0,\ y\geq0,\ 1\leq t=x+y\leq3\}.
\]
The pushforward of $dV_{\rm ext}$ is $(2\pi)^2\,dx\,dy$.
The extremal scalar curvature is affine on $P$, and the symmetry
interchanging $x$ and $y$ makes it a function $S=at+b$.
Let $d\sigma$ be lattice boundary measure. In the complex-scalar
normalization, the toric integration-by-parts identity for affine
functions $\ell$ reads
\[
 \int_P S\ell\,dx\,dy=\int_{\partial P}\ell\,d\sigma;
\]
compare \cite[Sections 5--7]{LeBrunToric}, whose Riemannian scalar
curvature is $2S$. The required moments are
\[
 |P|=4,\qquad \int_P t=\frac{26}{3},\qquad
 \int_P t^2=20,\qquad
 \int_{\partial P}1\,d\sigma=8,\qquad
 \int_{\partial P}t\,d\sigma=18.
\]
Taking $\ell=1,t$ therefore gives
\[
 \frac{26}{3}a+4b=8,\qquad
 20a+\frac{26}{3}b=18,
 \qquad
 S_{\rm ext}=\frac{6t+9}{11},\quad \bar S=2.
\]
The scalar curvature is nonconstant, and
\begin{equation}\label{ex:eq:f1-energy}
 \Ca(\omega_{\rm ext})
 =4\pi^2\int_1^3\left(\frac{6t-13}{11}\right)^2t\,dt
 =\frac{16\pi^2}{11}
 <\frac{64\pi^2}{3}=8\pi^2\delta(A).
\end{equation}
Thus the extremal energy has a strict margin below the long-time
threshold; the displayed integral is evaluated exactly.

We next choose a genuinely non-extremal invariant perturbation.
Let $Y=\nabla^{1,0}_{\omega_{\rm ext}}S_{\rm ext}$ and let
$\theta_Y(u)$ be its real Hamiltonian of zero mean with respect to
$dV_u$, for $\omega_u=\omega_{\rm ext}+\ddc u$.
For any $U(2)$-invariant extremal metric in this class, the extremal
field equals $Y$: the Futaki character and the nondegenerate pairing
of Hamiltonians determine it independently of the representative
\cite{FM}. Consequently an invariant metric is extremal precisely when
\[
 R(u)\coloneqq S_u-\bar S-\theta_Y(u)=0.
\]
At $u=0$, the variation of the Hamiltonian cancels the first-order
term in the scalar-curvature variation. More precisely,
$\theta_Y(u)=S_{\rm ext}-\bar S+Y(u)$ and
$\left.\frac{d}{d\varepsilon}\right|_0S_{\varepsilon\psi}
=-\mathcal D^*\mathcal D\psi+Y(\psi)$, giving
\[
 R(0)=0,\qquad
 DR_0(\psi)=-\mathcal D^*\mathcal D\psi,\qquad
 \mathcal D=\bar\partial\nabla^{1,0}_{\omega_{\rm ext}};
\]
see \cite[Sections 1 and 3]{HZFlow}. The kernel of this elliptic
operator is finite-dimensional, whereas the smooth invariant
potential space is infinite-dimensional. Choose a real invariant
$\psi$ outside that kernel. Then
$R(\varepsilon\psi)=-\varepsilon\mathcal D^*\mathcal D\psi
+O(\varepsilon^2)$ in $C^0$, so
$\omega_{\rm ext}+\varepsilon\ddc\psi$ is non-extremal for all
sufficiently small nonzero $\varepsilon$.
Positivity and the strict inequality in
\eqref{ex:eq:f1-energy} persist for small $\varepsilon$.
The resulting initial metric therefore satisfies
\eqref{ex:eq:general-startup}, and its flow is global with uniformly
bounded Sobolev constants.

Since $U(2)$ is a maximal compact subgroup of the reduced automorphism
group of $\mathbb F_1$, Huang--Zheng stability
\cite[Theorems 1.4--1.5]{HZFlow} applies: the modified Calabi flow
converges smoothly and exponentially to an extremal metric in the
same class. Fixed-class uniqueness identifies this metric with
$\omega_{\rm ext}$ up to a holomorphic automorphism, which
preserves the anticanonical polarization. The ordinary
Calabi flow consequently converges to its polarized isometry class,
in agreement with Theorem~\ref{app:thm:main}. This is convergence
after holomorphic pullback; a non-cscK extremal metric itself moves
under the ordinary flow by its extremal automorphisms.

These examples have smooth limits. They do not provide a flow from
a smooth polarized surface to a genuinely singular extremal orbifold.
Such an example would additionally require a fixed-class global flow
with a uniform Sobolev bound and identification of its singular limit;
an algebraic degeneration or a varying-class sequence alone would not
verify these requirements.

\vskip10pt
Haozhao Li, Institute of Geometry and Physics, and School of Mathematical
Sciences, University of Science and Technology of China, No. 96 Jinzhai
Road, Hefei, Anhui Province, 230026, China.
Email: \href{mailto:hzli@ustc.edu.cn}{hzli@ustc.edu.cn}\par\medskip

Bing Wang, Institute of Geometry and Physics, and School of Mathematical
Sciences, University of Science and Technology of China, Hefei 230026,
China, and Hefei National Laboratory, Hefei 230088, China.
Email: \href{mailto:topspin@ustc.edu.cn}{topspin@ustc.edu.cn}\par\medskip

Kai Zheng, University of Chinese Academy of Sciences, Beijing, 100190,
China.\\
Email: \href{mailto:KaiZheng@amss.ac.cn}{KaiZheng@amss.ac.cn}

\end{document}